\documentclass[10pt]{article}

\usepackage[all,cmtip]{xy}

\usepackage{setspace}
\usepackage{bbm,youngtab}
\usepackage{amscd,mathrsfs}
\usepackage{amssymb,amsmath,dsfont,amsfonts,amsthm}
\usepackage{enumerate}
\usepackage[mathscr]{eucal}
\usepackage{stmaryrd} 
\usepackage[cbgreek]{textgreek} 
\usepackage{booktabs}
\usepackage{multirow}
\usepackage{hhline} 

\usepackage{rotating}
\usepackage{ftnxtra}
\usepackage{fnpos}
\usepackage[titletoc,toc,title]{appendix}
\usepackage{graphicx}
\DeclareGraphicsExtensions{.pdf,.png,.jpg}
\usepackage{upgreek}
\usepackage{tikz-cd}

\usepackage{authblk}
\usepackage{bm}
\usepackage{caption}
\usepackage{mathtools}
\usepackage{soul}            
\usepackage{xcolor}
\usepackage[colorlinks,citecolor=blue,linkcolor=blue]{hyperref}

\usepackage{natbib}

\definecolor{myColor}{RGB}{5,175,255}

\numberwithin{equation}{section}

\theoremstyle{definition}
\newtheorem{remark}{Remark}

\theoremstyle{definition}

\usepackage{enumitem}

\title{
\vspace{-3cm}
\textbf{First-Order Compatible-Strain Mixed Quadrilateral Finite Elements for $2$D Nonlinear Elasticity}}

\author[1,2]{Mohsen Jahanshahi\thanks{Corresponding author, e-mail: mohsen.jahanshahi@mcgill.ca}}

\author[1]{Damiano Pasini}
\affil[1] {\small \textit{Department of Mechanical Engineering, McGill University, Montreal, Quebec H3A 0C3, Canada}}
\affil[2]{\small \textit{Department of Mechanical, Environmental and Civil Engineering, Mayfield College of Engineering, Tarleton State University, Stephenville, TX 76401, USA}}

\author[3,4]{Arash Yavari}
\affil[3]{\small \textit{School of Civil and Environmental Engineering, Georgia Institute of Technology, Atlanta, GA 30332, USA}}
\affil[4]{\small \textit{The George W. Woodruff School of Mechanical Engineering, Georgia Institute of Technology, Atlanta, GA 30332, USA}}

\begin{document}

\maketitle

\begin{abstract}
Compatible-strain mixed finite elements (CSMFEs) use the differential complex of nonlinear elasticity to construct discretizations that preserve the underlying topological structure. Existing developments of CSMFEs have focused on simplicial meshes, for which efficient formulations have been developed for first- and second-order elements and successfully applied to compressible and incompressible nonlinear elasticity problems. In this paper, we develop two compatible-strain mixed formulations for quadrilateral elements applicable to compressible and incompressible nonlinear elasticity. For quadrilateral elements, the Piola transformation preserves vector fields tangent and normal to element edges only for special classes of elements, such as rectangles, parallelograms, and trapezoids. Consequently, it cannot be used to construct compatible shape functions for general quadrilateral elements. To overcome this limitation, the shape functions are computed directly in the physical space using numerical integration. N\'ed\'elec shape functions of the first kind are used to construct the interpolation of the displacement gradient. We also introduce a new class of shape functions for the stress tensor that is compatible with the displacement and displacement-gradient discretizations. The formulation for compressible elasticity follows the framework previously developed for first- and second-order simplicial CSMFEs, whereas a new formulation is proposed for incompressible elasticity. The incompressible formulation employs element-level condensation of the pressure field and therefore does not increase the number of global degrees of freedom relative to the corresponding compressible formulation. These developments extend the compatible-strain mixed finite element framework from simplicial to general quadrilateral meshes for both compressible and incompressible nonlinear elasticity. Numerical examples demonstrate that the proposed quadrilateral elements can successfully solve problems that require second-order simplicial CSMFEs while using fewer degrees of freedom.
\end{abstract}

\begin{description}
\item[Keywords:] Mixed finite elements, finite element exterior calculus, compatible-strain finite elements, nonlinear elasticity, N{\'e}d{\'e}lec shape functions.
\end{description}

\tableofcontents

\section{Introduction}
\label{sec:1}

The displacement-based Finite Element Method (FEM) is widely used in engineering analyses and simulations. Its popularity stems from the simplicity of its formulation, as discussed in detail in \citep{BAT96,ZIE05,HUG00,BEL05,WRI08}, as well as its computational efficiency and modest computational requirements. However, the stress tensor and displacement gradient are generally discontinuous across element boundaries in displacement-based finite elements. These formulations also perform poorly when modeling incompressible or nearly incompressible solids \citep{GLO84,GLO88,SIM91}. Mixed finite elements were introduced in order to address these limitations. The formulation of mixed finite elements for linear elasticity is discussed in detail in \citep{WAS75,BRA07,BRE07}. Early developments of mixed finite elements for linear elasticity include the pioneering works of \citet{PIA64} and \citet{PIA84}. For finite-strain applications, mixed formulations were developed by introducing a multiplicative decomposition of the deformation gradient into volumetric and isochoric parts and treating displacement, pressure, and dilation as independent fields \citep{SIM85,SIM88a,SIM88b,SIM91}. Enhanced strain methods and methods based on incompatible modes were subsequently developed within the mixed formulation framework \citep{SIM90,SIM92,SIM93}. Various stabilization techniques were later proposed to address the hourglass instabilities observed in first-order mixed finite elements \citep{KOR96,WRI96,GLA97,REE00}. \citet{SCH11} developed a mixed finite element formulation for quasi-incompressible finite elasticity based on approximations of the minors of the deformation gradient. This formulation was later extended by treating the deformation gradient, its adjoint, and its determinant as independent kinematic variables \citep{BON15}. \citet{NEU21} introduced three mixed formulations for nonlinear elasticity based on different combinations of displacement, deformation gradient, right Cauchy-Green strain, and first and second Piola-Kirchhoff stress tensors as independent fields. In the context of finite-strain elastoplasticity, \citet{SIM88a,SIM88b} utilized the multiplicative decomposition of the deformation gradient into elastic and plastic parts to formulate a mixed variational framework. \citet{JAH15} used a three-field Hu-Washizu principle together with a similar decomposition to develop an integration algorithm for $J_{2}$ plasticity. \citet{KHO16} extended classical plasticity to nanostructures using a mixed formulation and the reversed multiplicative decomposition of the deformation gradient into plastic and elastic parts.

A distinct class of mixed finite elements for linear elasticity is based on the differential complexes of linear elasticity \citep{ARN02,ARN06,ARN18}. In particular, \citet{ARN06} proposed a discretization of the linear elasticity complex that preserves its underlying topological structure. The Finite Element Exterior Calculus (FEEC) developed in \citep{ARN06,ARN10,ARN18} uses tools from geometry and topology to construct numerical methods for classes of partial differential equations that are often challenging for conventional finite element techniques. For nonlinear elasticity applications, \citet{ANG15,ANG16} introduced a differential complex suitable for describing the kinematics and kinetics of large deformations. Similar discrete complexes can be constructed using vector-valued shape functions introduced by \citet{RAV77}, \citet{BRE85}, and \citet{NED80,NED86}. Building on this framework, \citet{ANG17} used displacement, displacement gradient, and the first Piola-Kirchhoff stress tensor as independent fields in a Hu-Washizu functional to develop a mixed formulation for two-dimensional compressible nonlinear elasticity. Because the displacement gradient satisfies the classical Hadamard jump condition for compatibility of non-smooth strain fields, they termed their method the \textit{Compatible-Strain Mixed Finite Element (CSMFE)}. \citet{DHA22a} used Cartan's moving frames to formulate a mixed variational approach to nonlinear elasticity and implemented it within the FEEC framework. This formulation was subsequently extended to three-dimensional nonlinear elasticity \citep{DHA22b}. \citet{SHO18} used a four-field Hu-Washizu functional involving displacement, displacement gradient, the first Piola-Kirchhoff stress tensor, and pressure to extend CSMFEs to incompressible nonlinear elasticity. Later, \citet{SHO19} extended this formulation to three-dimensional compressible and incompressible nonlinear elasticity. \citet{JAH22} developed a compatible mixed finite element formulation for first-order simplex elements using mid-nodes to enforce continuity constraints across element boundaries. This formulation was subsequently extended to second-order simplex elements by \citet{JAH25}. It was shown that the resulting first- and second-order CSMFEs possess excellent computational efficiency while effectively addressing difficulties commonly encountered in conventional first-order finite elements. Their performance was also compared with that of conventional first- and second-order finite elements and mixed displacement-pressure $U/P$ formulations.

In this paper, we formulate two mixed formulations for the analysis of two-dimensional compressible and incompressible solids undergoing finite deformations. The formulation of CSMFEs, originally developed for simplicial elements, is extended to quadrilateral elements for nonlinear elasticity applications. It is important to note that for existing $h\left(\operatorname{div}\right)$ or $h\left(\operatorname{curl}\right)$ finite elements either the normal component of the stress tensor is single-valued over an edge shared by two elements or the tangent component of the displacement gradient is single-valued over this edge \citep{ARN05,FAL11,ROG10,SCH22}. This is in contrast to CSMFEs where both the normal component of the stress tensor and the tangent component of the displacement gradient are single-valued over the common edge. Following the approach introduced in \citep{JAH22,JAH25}, the first formulation employs displacement, displacement gradient, and stress tensor as independent fields in a Hu-Washizu functional for compressible nonlinear elasticity. The second formulation is based on a five-field functional involving displacement, displacement gradient, stress tensor, dilation, and pressure, and is designed for incompressible nonlinear elasticity. For general quadrilateral elements in the physical space, the covariant and contravariant Piola transformations \citep{ROG10,AZN22} do not preserve the property that vector fields tangent and normal to a given edge in the natural coordinate system remain tangent and normal to the corresponding edge after transformation. Consequently, these transformations cannot be used to construct compatible shape functions for arbitrary quadrilateral elements. To overcome this difficulty, the shape functions used to interpolate the displacement gradient and stress tensor are computed directly in the physical space using numerical integration. The N\'ed\'elec shape functions of the first kind \citep{NED80,NED86} are used to construct the displacement-gradient interpolation. We also introduce a new class of shape functions for the stress tensor that is compatible with the degrees of freedom associated with the displacement and displacement-gradient fields. Finally, sufficient and necessary conditions for the invertibility of the global stiffness matrix are derived and analyzed in detail.

\paragraph{Contributions of this paper.}
The main contributions of this paper are as follows.

\begin{itemize}[topsep=2pt,noitemsep,leftmargin=10pt]

\item The CSFME formulation, originally developed for simplicial elements, is extended to quadrilateral elements. Numerical experiments indicate that certain problems cannot be solved over the entire range of applied loads using first-order simplicial elements with $24$ degrees of freedom \citep{JAH22}. In such cases, second-order simplicial elements with $54$ degrees of freedom are required \citep{JAH25}. The quadrilateral element developed in this work, with $48$ degrees of freedom, successfully solves these problems without exhibiting the deficiencies associated with first-order simplicial elements. These results suggest that the proposed quadrilateral formulation provides a more efficient alternative to second-order simplicial elements.

\item A novel formulation for modeling incompressible solids is presented. Unlike conventional mixed formulations, it does not increase the number of degrees of freedom of the element, as the pressure degree of freedom is condensed out at the element level. Consequently, incompressible solids can be modeled with only a minimal computational overhead compared to compressible solids.

\item The Piola transformation applied to quadrilateral elements preserves the tangential or normal character of the corresponding vector field along an edge only for certain classes of elements, such as rectangles, parallelograms, and trapezoids. To circumvent this limitation, the shape functions for the displacement gradient and stress tensor are constructed directly in the physical space using numerical integration. Consequently, the proposed formulation can be applied to general quadrilateral meshes, including irregular meshes.

\item Conditions that guarantee the invertibility of the overall stiffness matrix are identified in terms of the constituent stiffness matrices of the element. The shape functions used to interpolate the stress tensor are then constructed to satisfy these conditions, thereby ensuring the invertibility of the resulting element stiffness matrix.

\end{itemize}

This paper is organized as follows. In \S\ref{sec:2}, we introduce the shape functions used to interpolate the displacement gradient and stress tensor. In \S\ref{sec:3}, we discuss the degrees of freedom associated with the displacement, displacement gradient, and stress tensor fields, and present the interpolation procedures for the displacement gradient and stress tensor. This section also presents the formulations of the compressible and incompressible CSMFEs and concludes with a discussion of the conditions required for solvability of the resulting system of equations. Several numerical examples are presented in \S\ref{sec:4}. The performance of the proposed mixed formulations is compared with that of existing first-order mixed finite elements. In particular, to assess the convergence properties of the incompressible formulation, contours of various stress components and pressure are compared with those obtained using a second-order mixed $U/P$ finite element. Finally, concluding remarks are given in \S\ref{sec:5}.

\section{Vector-Valued Polynomial Shape Functions}
\label{sec:2}

In this section, we discuss in detail the vector-valued shape functions used to interpolate the displacement gradient and stress tensors. To achieve the optimal convergence rate for rectangular meshes with bilinear displacement shape functions, the N\'ed\'elec shape functions of the first kind are used for the interpolation of the displacement gradient. Since the Piola transformation does not preserve the property that vector-valued shape functions remain tangent or normal to a given edge of a general element in the physical space, the shape functions are computed directly in the physical space. The procedure for computing the coefficients of the local shape functions on individual elements, as well as the construction of global shape functions on edges with prescribed global orientations, is described. Next, we introduce a new class of shape functions for interpolating the stress tensor. The degree of these shape functions is chosen so that they are compatible with the degrees of freedom associated with the displacement and displacement-gradient fields. In other words, these shape functions are constructed so that the resulting stiffness matrices satisfy the inf--sup, or Ladyzhenskaya--Babu\v{s}ka--Brezzi, conditions \citep{BOF13}. This point is discussed further in \S\ref{sec:3}. The construction of global shape functions for interpolating the stress tensor is also presented.

\begin{figure}
\centering
\includegraphics[width=0.5\textwidth]{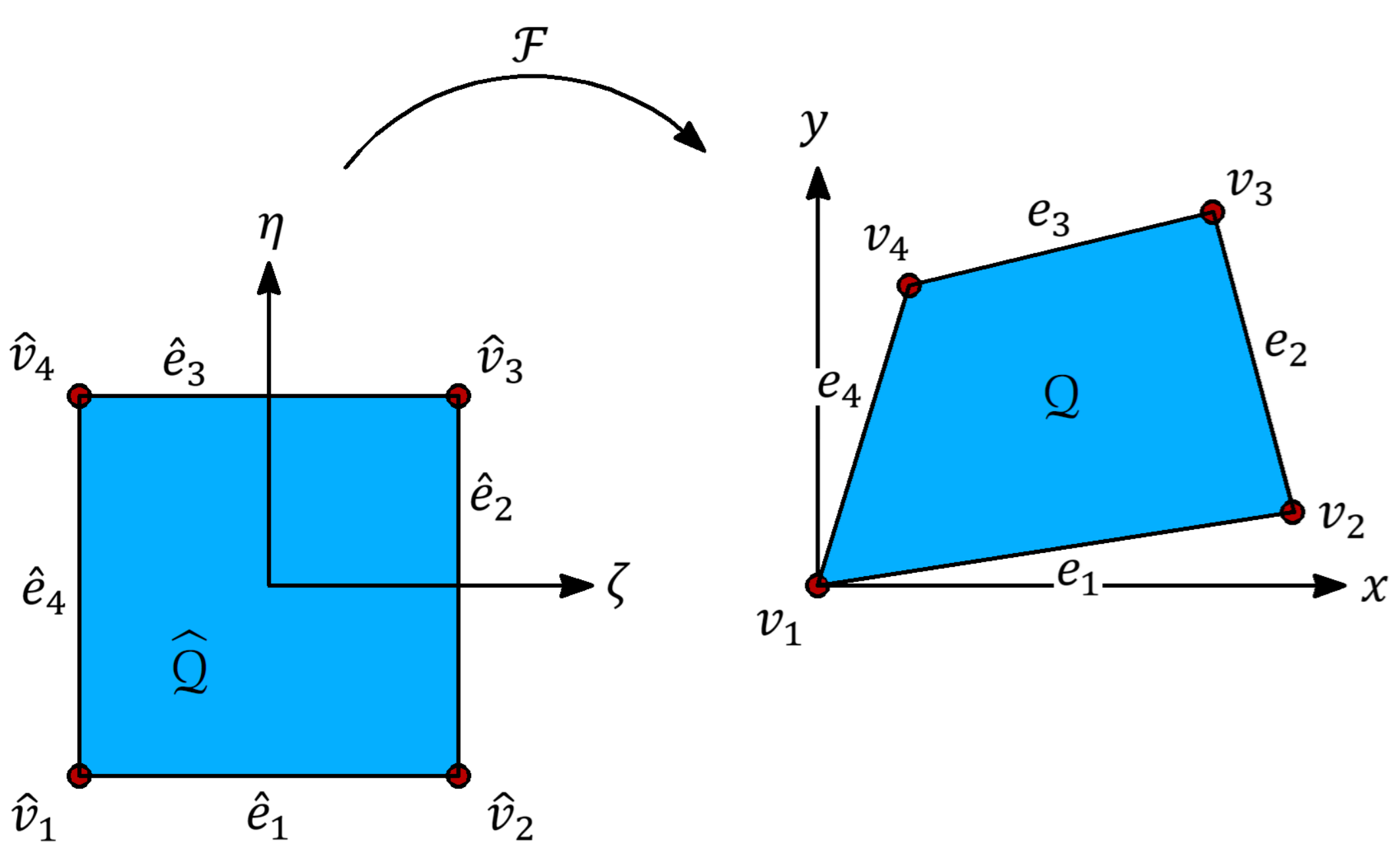}
\vskip 0.1in
\caption{A square element $\widehat{\mathscr{Q}}$ with the vertices and edges $\hat{v}_{i},\hat{e}_{i},i=1,\dots,4$ in the natural coordinates $\zeta$ and $\eta$ is mapped via $\mathcal{F}$ to the element $\mathscr{Q}$ with the vertices and edges $v_{i},e_{i},~i=1,\dots,4$ in the physical space.}
\label{fig:21}
\end{figure}

In the following discussion, the domain $\mathcal{B}$ with boundary $\partial\mathcal{B}$ is assumed to be discretized into a collection $\mathcal{Q}$ of convex quadrilateral elements $\mathscr{Q}$ in the physical space. Let $\mathcal{E}$ denote the set of all edges of the mesh and let $\mathcal{E}^{\partial}$ denote the set of edges lying on the boundary $\partial\mathcal{B}$. The set of interior edges is then given by $\mathcal{E}^{i}=\mathcal{E}\setminus\mathcal{E}^{\partial}$. As illustrated in Figure~\ref{fig:21}, the mapping $\mathcal{F}:[-1,1]^{2}\rightarrow\mathbb{R}^{2}$ maps the reference square element $\widehat{\mathscr{Q}}$ in the natural coordinates $(\zeta,\eta)$ to the quadrilateral element $\mathscr{Q}$ in the physical space. The vertices and edges of the reference element are denoted by $\hat{v}_{i}$ and $\hat{e}_{i}$, $i=1,\ldots,4$, whereas $v_{i}$ and $e_{i}$, $i=1,\ldots,4$, denote the corresponding vertices and edges of $\mathscr{Q}$ in the physical space. The vertices and edges of both elements are oriented counterclockwise. We further assume that a unique global orientation is assigned to each edge $e\in\mathcal{E}$. The vector-valued shape functions used to interpolate the displacement gradient are denoted by $\mathbf{v}\in\mathcal{P}_{3}(T\mathbb{R}^{2})$, while those used to interpolate the stress tensor are denoted by $\mathbf{w}\in\mathcal{P}_{2}(T\mathbb{R}^{2})$, where $\mathcal{P}_{3}(T\mathbb{R}^{2})$ and $\mathcal{P}_{2}(T\mathbb{R}^{2})$ are the spaces of polynomial vector fields of degree at most $3$ and $2$, respectively. The jump of the tangential component of $\mathbf{v}$ across an edge $e$ shared by two quadrilateral elements $\mathscr{Q}_{1},\mathscr{Q}_{2}\in\mathcal{Q}$ is defined as \citep{ANG17,JAH22,JAH25}:
\begin{equation}\label{eqn:2.1}
{\llbracket\mathrm{t}\mathbf{v}\rrbracket}_{e}=
{\mathbf{v}|}_{\mathscr{Q}_{1}}\cdot\mathbf{t}^{e}-
{\mathbf{v}|}_{\mathscr{Q}_{2}}\cdot\mathbf{t}^{e}\,,
\end{equation}
where $\mathbf{t}^{e}$ is the global unit vector tangent to the edge $e$. Similarly, denoting $\mathbf{n}^{e}$ as the global unit vector normal to the edge, the jump of the normal component of $\mathbf{w}$ across $e$ can be defined as:
\begin{equation}\label{eqn:2.2}
{\llbracket\mathrm{n}\mathbf{w}\rrbracket}_{e}=
{\mathbf{v}|}_{\mathscr{Q}_{1}}\cdot\mathbf{n}^{e}-
{\mathbf{v}|}_{\mathscr{Q}_{2}}\cdot\mathbf{n}^{e}\,.
\end{equation}
Using equations \eqref{eqn:2.1} and \eqref{eqn:2.2}, the spaces of polynomial vector fields of degree at most $r$ with zero jumps for tangent and normal components across an internal edge $e\in\mathcal{E}^{i}$ are defined, respectively, as
\begin{equation}\label{eqn:2.3}
\begin{aligned}
\mathcal{P}^{c}_{r}\left(T\mathcal{Q}\right)
&=\left\{\mathbf{v}\in\mathcal{P}_{r}\left(T\mathcal{Q}\right)\Big|
{\llbracket\mathrm{t}\mathbf{v}\rrbracket}_{e}=0\,,\,\forall e\in\mathcal{E}^{i}\right\}\,,\\
\mathcal{P}^{d}_{r}\left(T\mathcal{Q}\right)
&=\left\{\mathbf{w}\in\mathcal{P}_{r}\left(T\mathcal{Q}\right)\Big|
{\llbracket\mathrm{n}\mathbf{w}\rrbracket}_{e}=0,\,\forall e\in\mathcal{E}^{i}\right\}\,.
\end{aligned}
\end{equation}

\subsection{Space of c-conformal shape functions}
\label{sec:2.1}

The N\'ed\'elec shape functions of the first kind $\mathbf{v}\in Q_{l-1,l}\times Q_{l,l-1}\subset\mathcal{P}_{2l-1}\left(T\mathbb{R}^2\right)$ are used to interpolate the displacement gradient. Here, $Q_{r,s}\subset\mathcal{P}_{r+s}\left(\mathbb{R}\right)$ denotes the space of polynomials of the form $x^{\alpha}y^{\beta}$, $\alpha=0,\dots,r$, $\beta=0,\dots,s$ \citep{NED80,NED86}. Following \citet{JAH25}, we consider the case $l=2$. The corresponding local shape functions $\mathbf{v}$ associated with the edges $e_i$, $i=1,\dots,4$, are given by:
\begin{equation}\label{eqn:2.4}
\mathbf{v}^{\mathscr{Q},e_{i}}_{j}=
\begin{bmatrix}
  a^{j}_{1}+b^{j}_{1}x+c^{j}_{1}y+d^{j}_{1}xy+e^{j}_{1}y^{2}+f^{j}_{1}xy^{2} \\
  a^{j}_{2}+b^{j}_{2}y+c^{j}_{2}x+d^{j}_{2}xy+e^{j}_{2}x^{2}+f^{j}_{2}x^{2}y
\end{bmatrix}\,,
\end{equation}
where $j=1,2$ refers to the index of the shape function on the edge $e_{i}$. To complete the space of shape functions that interpolate the displacement gradient, it is also necessary to define the following shape functions on the quadrilateral element $\mathscr{Q}\in\mathcal{Q}$:
\begin{equation}\label{eqn:2.5}
\mathbf{v}^{\mathscr{Q}}_{k}=
\begin{bmatrix}
  a^{k}_{1}+b^{k}_{1}x+c^{k}_{1}y+d^{k}_{1}xy+e^{k}_{1}y^{2}+f^{k}_{1}xy^{2} \\
  a^{k}_{2}+b^{k}_{2}y+c^{k}_{2}y+d^{k}_{2}xy+e^{k}_{2}x^{2}+f^{k}_{2}x^{2}y
\end{bmatrix}\,,
\end{equation}
where $k=1,\dots,4$ denotes the index of the shape function on $\mathscr{Q}$. From equations \eqref{eqn:2.4} and \eqref{eqn:2.5}, it follows that $\mathbf{v}^{\mathscr{Q},e_{i}}_{j}, \mathbf{v}^{\mathscr{Q}}_{k}\in\mathcal{P}_{3}\left(T\mathbb{R}^{2}\right)$. To determine the coefficients of the shape functions $\mathbf{v}^{\mathscr{Q},e_{i}}_{j}$, $i=1,\dots,4$, $j=1,2$, and $\mathbf{v}^{\mathscr{Q}}_{k}$, $k=1,\dots,4$, we define the following degrees of freedom on the edges $e_{i}$, $i=1,\dots,4$, and on the quadrilateral $\mathscr{Q}$ itself \citep{JAH25}:
\begin{equation}\label{eqn:2.6}
\begin{aligned}
\phi^{\mathscr{Q},e_{i}}_{j}\left(\mathbf{v}\right)&=\int_{e_{i}}\left(\mathbf{v}\cdot\mathbf{t}_{i}\right)h_{j}ds=
\int_{\hat{e}_{i}}\big[\mathbf{v}\left(\mathcal{F}\left(\hat{\mathbf{x}}\right)\right)\cdot\mathbf{t}_{i}\big]\hat{h}_{j}
\left(\hat{\mathbf{x}}\right)\left|\mathcal{F}|_{\hat{e}_{i}}\right|d\hat{s}
\,,&& i=1,\dots,4,j=1,2\,, \\
\phi^{\mathscr{Q}}_{k}\left(\mathbf{v}\right)&=\int_{\mathscr{Q}}\left(\mathbf{v}\cdot\mathbf{b}^{\mathscr{Q}}_{k}\right)dA=
\int_{\widehat{\mathscr{Q}}}\big[\mathbf{v}\left(\mathcal{F}\left(\hat{\mathbf{x}}\right)\right)\cdot
\mathbf{b}^{\mathscr{Q}}_{k}\left(\mathcal{F}\left(\hat{\mathbf{x}}\right)\right)\big]
\big|\mathcal{F}\left(\hat{\mathbf{x}}\right)\big|d\widehat{A}\,,&& k=1,\dots,4\,,
\end{aligned}
\end{equation}
where $\mathbf{t}_{i}$ is the unit vector tangent to the edge $e_{i}$, $\mathcal{F}$ is the map from the reference square in the natural coordinate system to the quadrilateral element $\mathscr{Q}$ in the physical space (see Figure~\ref{fig:21}) and $\mathcal{F}|_{\hat{e}_{i}}$ is its restriction to the edge $\hat{e}_{i}$. Furthermore, $\left|\mathcal{F}|_{\hat{e}_{i}}\right|$ and $\big|\mathcal{F}\left(\hat{\mathbf{x}}\right)\big|$ are the Jacobian of the transformations from the local coordinate system to the physical space. Following \citet{SCH22}, the basis functions $\hat{h}_{j},j=1,2$ on the edge $\hat{e}_{i}$ are considered as follows:
\begin{equation}\label{eqn:2.7}
\hat{h}_{1}=\frac{1}{2}\left(1-s\right),\qquad\hat{h}_{2}=\frac{1}{2}\left(1+s\right)\,.
\end{equation}
We use the following base vectors to define the degrees of freedom $\phi^{\mathscr{Q}}_{k},k=1,\dots,4$ in \eqref{eqn:2.6}$_{2}$:
\begin{equation}\label{eqn:2.8}
\mathbf{b}^{\mathscr{Q}}_{1}=
\begin{bmatrix}
  \frac{1}{2}\left(1-x\right) \\
  0
\end{bmatrix},\qquad
\mathbf{b}^{\mathscr{Q}}_{2}=
\begin{bmatrix}
  \frac{1}{2}\left(1+x\right) \\
  0
\end{bmatrix},\qquad
\mathbf{b}^{\mathscr{Q}}_{3}=
\begin{bmatrix}
  0 \\
  \frac{1}{2}\left(1-y\right)
\end{bmatrix},\qquad
\mathbf{b}^{\mathscr{Q}}_{4}=
\begin{bmatrix}
  0 \\
  \frac{1}{2}\left(1+y\right)
\end{bmatrix}\,.
\end{equation}
The coefficients of the shape functions $\mathbf{v}^{\mathscr{Q},e_{i}}_{j}$ are obtained by requiring the following conditions to be satisfied \citep{JAH25}:
\begin{equation}\label{eqn:2.9}
\phi^{\mathscr{Q},e_{l}}_{m}\left(\mathbf{v}^{\mathscr{Q},e_{i}}_{j}\right)=
\begin{cases}
  1,&\text{if } l=i,~m=j\,, \\
  0,&\text{otherwise}\,,
\end{cases}
\quad\text{and}\qquad
\phi^{\mathscr{Q}}_{n}\left(\mathbf{v}^{\mathscr{Q},e_{i}}_{j}\right)=0\,.
\end{equation}
Similarly, the coefficients of the shape functions $\mathbf{v}^{\mathscr{Q}}_{k}$ can be determined by imposing the following conditions:
\begin{equation}\label{eqn:2.10}
\phi^{\mathscr{Q},e_{l}}_{m}\left(\mathbf{v}^{\mathscr{Q}}_{k}\right)=0,\qquad
\phi^{\mathscr{Q}}_{n}\left(\mathbf{v}^{\mathscr{Q}}_{k}\right)=
\begin{cases}
  1,&\text{if } n=k\,,\\
  0,&\text{otherwise}\,.
\end{cases}
\end{equation}

\begin{remark}\label{rem:2.1}
The coefficients of the shape functions are computed by solving matrix equations for the unknown coefficients $a^{j}_{1},\dots,f^{j}_{2}$ or $a^{k}_{1},\dots,f^{k}_{2}$. The coefficient matrix is obtained by substituting \eqref{eqn:2.4} or \eqref{eqn:2.5} into the right-hand side of \eqref{eqn:2.6} and numerically integrating the monomials $1,x,y,\dots,x^{2}y$ that appear in $\mathbf{v}$. For all $12$ matrix equations ($8$ for $\mathbf{v}^{\mathscr{Q},e_{i}}_{j}$ and $4$ for $\mathbf{v}^{\mathscr{Q}}_{k}$), the coefficient matrix is identical. The right-hand side of each equation is determined by the conditions in \eqref{eqn:2.9} and \eqref{eqn:2.10}.
\end{remark}

\begin{figure}
\centering
\includegraphics[width=0.75\textwidth]{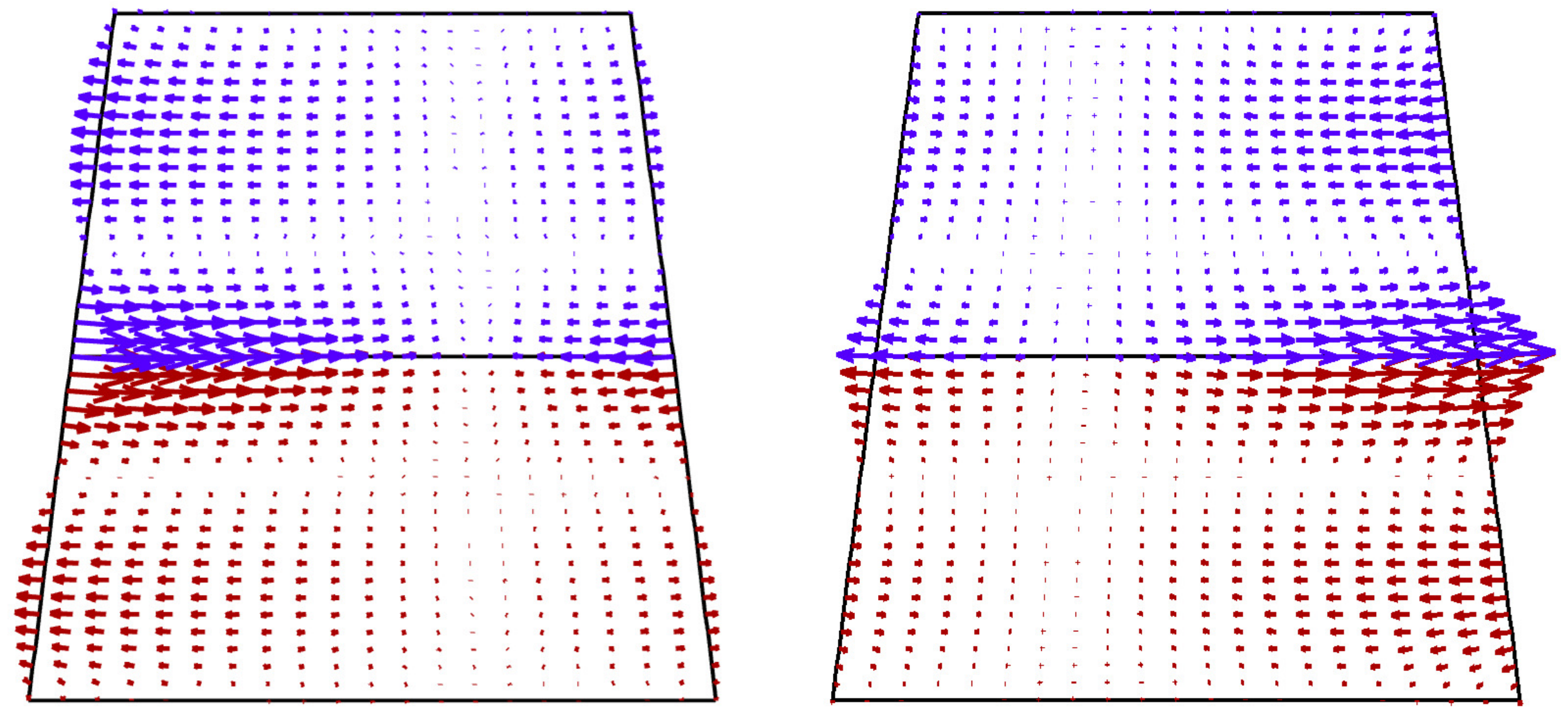}
\vskip 0.1in
\caption{Global shape functions $\mathbf{V}^{e}_{j},j=1,2$ on an edge $e$ shared by two quadrilaterals $\mathscr{Q}_{1}$ and $\mathscr{Q}_{2}$. It is observed that the tangent component of the vector fields in the two quadrilaterals is single-valued over the common edge $e$.}
\label{fig:22}
\end{figure}

\begin{figure}
\centering
\includegraphics[width=0.75\textwidth]{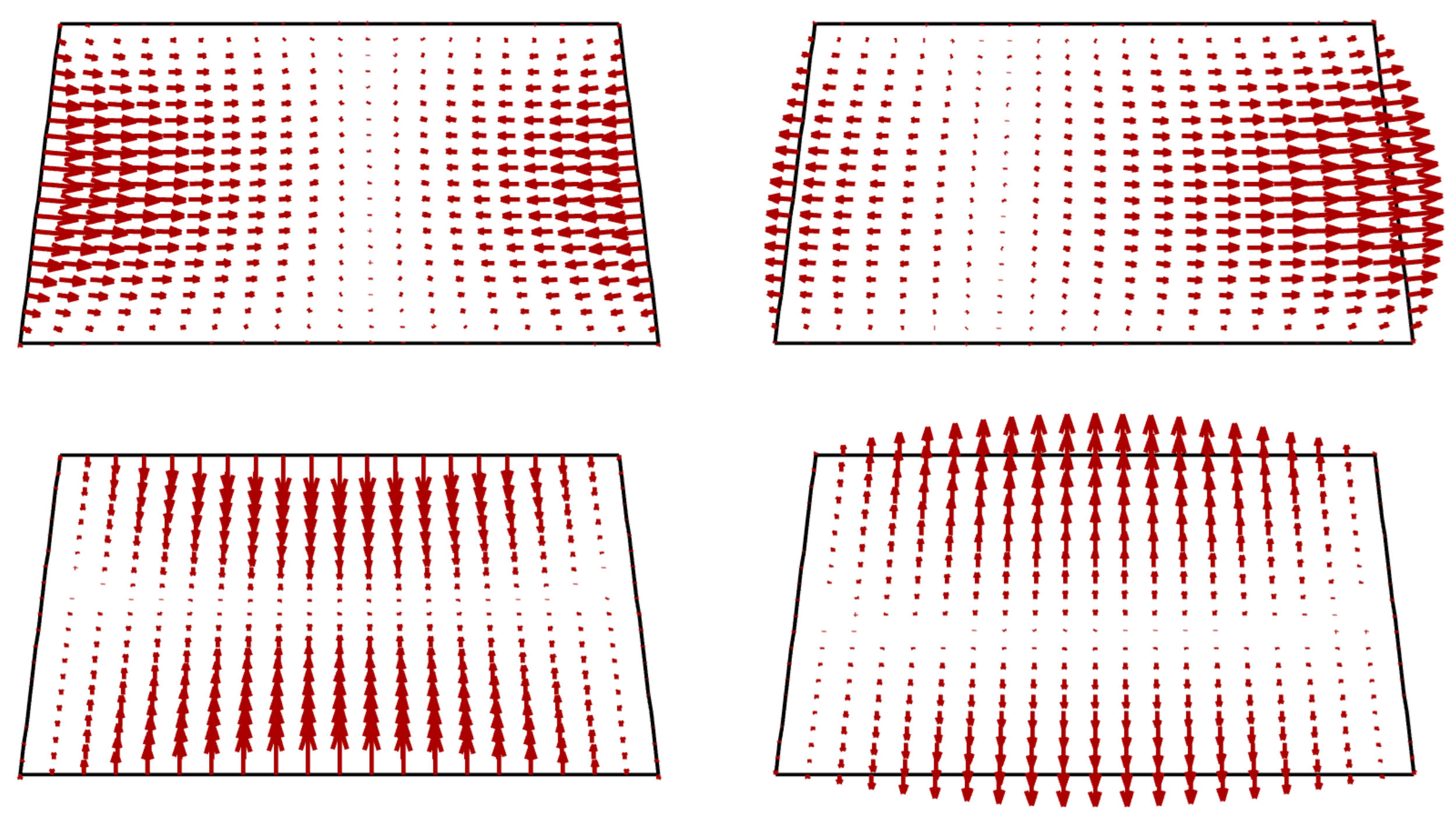}
\vskip 0.1in
\caption{Global shape functions $\mathbf{V}^{\mathscr{Q}}_{k},k=1,\dots,4$ on the quadrilateral $\mathscr{Q}$. We note that the vector field represented by the shape functions either vanishes on the edges of the element or it is normal to them.}
\label{fig:23}
\end{figure}

To enforce the jump condition for the tangential component of $\mathbf{v}^{\mathscr{Q},e_{i}}_{j}$ across an edge $e$ shared by two quadrilaterals $\mathscr{Q}_{1}$ and $\mathscr{Q}_{2}$, namely the condition in \eqref{eqn:2.3}$_{1}$, we define the global shape functions $\mathbf{V}^{e}_{j}$, $j=1,2$, as follows:
\begin{equation}\label{eqn:2.11}
\mathbf{V}^{e}_{j}|_{\mathscr{Q}_{l}}=
\begin{cases}
  c^{l}_{i}\mathbf{v}^{\mathscr{Q}_{l},e_{i}}_{j}\,,&
  \text{if }\mathscr{Q}_{l}=\mathscr{Q}_{1},\mathscr{Q}_{2}\,,\\
  0,&\text{otherwise}\,,
\end{cases}
\end{equation}
where $e_{i}\in\partial\mathscr{Q}_{l}$ is the same edge with the global identifier $e$ and,
\begin{equation}\label{eqn:2.12}
c^{l}_{i}=\mathbf{t}^{e}\cdot\mathbf{t}^{l}_{i}\,.
\end{equation}
In \eqref{eqn:2.12}, $\mathbf{t}^{e}$ denotes the global unit tangent vector associated with the edge $e$, while $\mathbf{t}^{l}_{i}$ denotes the local unit tangent vector associated with the edge $e_{i}\in\partial\mathscr{Q}_{l}$. It should be noted that the global shape functions $\mathbf{V}^{\mathscr{Q}}_{k}$, $k=1,\dots,4$, associated with the quadrilateral $\mathscr{Q}$ always satisfy:
\begin{equation}\label{eqn:2.13}
\mathbf{V}^{\mathscr{Q}}_{k}=\mathbf{v}^{\mathscr{Q}}_{k},k=1,\dots,4\,.
\end{equation}
The global shape functions $\mathbf{V}^{e}_{j}$, $j=1,2$, associated with an edge $e$ shared by two quadrilaterals $\mathscr{Q}_{1}$ and $\mathscr{Q}_{2}$ are illustrated in Figure~\ref{fig:22}. As seen in the figure, the tangential component of the corresponding vector field is single-valued across the common edge $e$. Figure~\ref{fig:23} shows the global shape functions $\mathbf{V}^{\mathscr{Q}}_{k}$, $k=1,\dots,4$, associated with the quadrilateral $\mathscr{Q}$. For these shape functions, the corresponding vector field either vanishes on the edges of the element or is normal to them.


\subsection{Space of d-conformal shape functions}
\label{sec:2.2}

\begin{figure}
\centering
\includegraphics[width=0.75\textwidth]{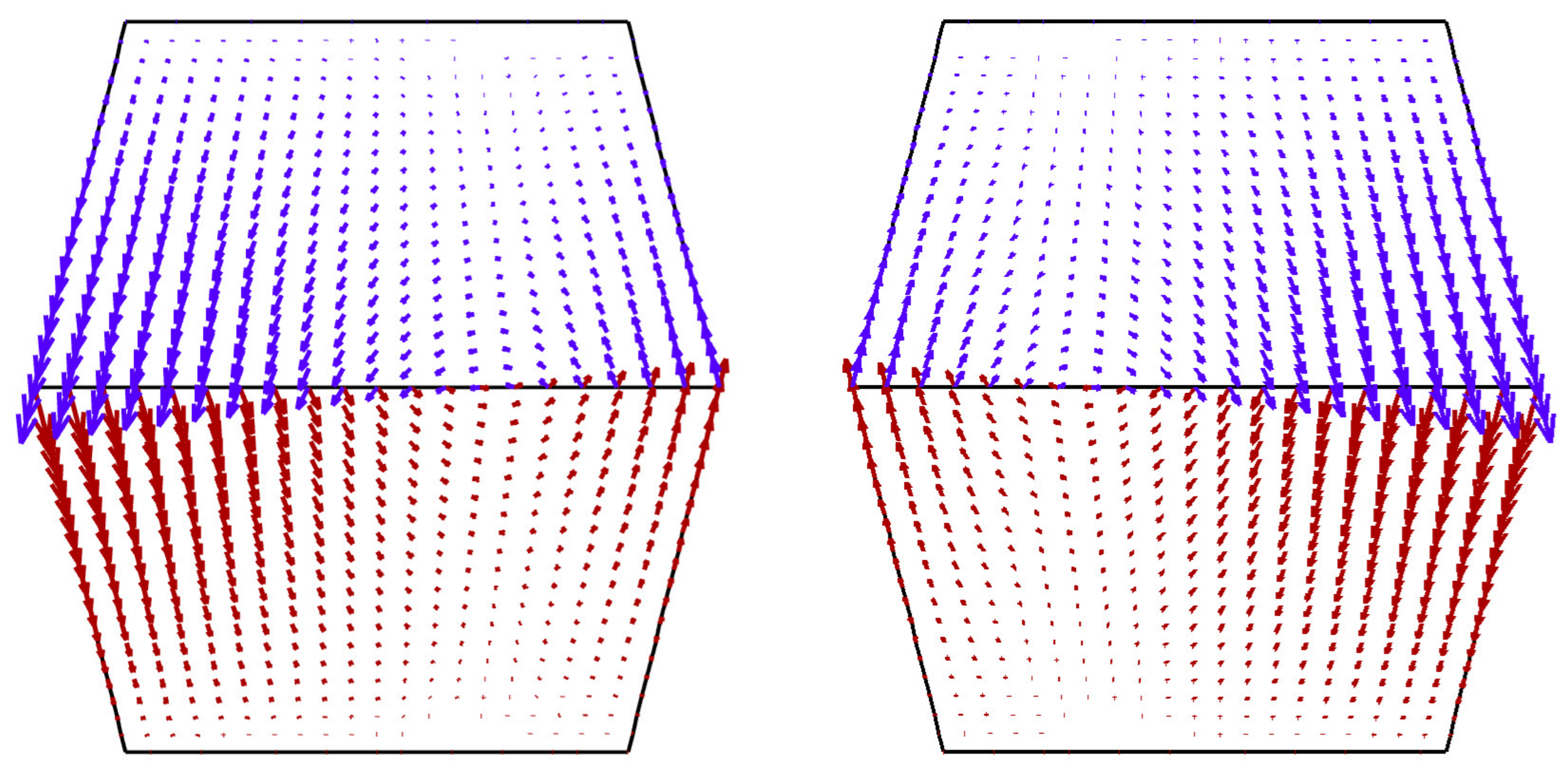}
\vskip 0.1in
\caption{Global shape functions $\mathbf{W}^{e}_{j},j=1,2$ on an edge $e$ shared by two quadrilaterals $\mathscr{Q}_{1}$ and $\mathscr{Q}_{2}$. It is observed that the normal component of the vector fields in the two quadrilaterals is single-valued over the common edge $e$.}
\label{fig:24}
\end{figure}

To interpolate the stress tensor, we define the local shape functions $\mathbf{w}^{\mathscr{Q},e_{i}}_{j}\in\mathcal{P}_{2}\left(T\mathbb{R}^{2}\right)$ on the edges $e_{i},i=1,\dots,4$ as follows:
\begin{equation}\label{eqn:2.14}
\mathbf{w}^{\mathscr{Q},e_{i}}_{j}=
\begin{bmatrix}
  a^{j}_{1}+b^{j}_{1}x+c^{j}_{1}y+d^{j}_{1}xy \\
  a^{j}_{2}+b^{j}_{2}x+c^{j}_{2}y+d^{j}_{2}xy
\end{bmatrix}\,,
\end{equation}
where $j=1,2$ refers to the index of the shape function on the edge $e_{i}$. To determine the coefficients of the shape functions $\mathbf{w}^{\mathscr{Q},e_{i}}_{j},i=1,\dots,4,j=1,2$, it is necessary to define the following degrees of freedom on the edges $e_{i},i=1,\dots,4$:
\begin{equation}\label{eqn:2.15}
\varphi^{\mathscr{Q},e_{i}}_{j}\left(\mathbf{w}\right)=\int_{e_{i}}\left(\mathbf{w}\cdot\mathbf{n}_{i}\right)h_{j}ds=
\int_{\hat{e}_{i}}\big[\mathbf{w}\left(\mathcal{F}\left(\hat{\mathbf{x}}\right)\right)\cdot\mathbf{n}_{i}\big]\hat{h}_{j}
\left(\hat{\mathbf{x}}\right)\big|\mathcal{F}|_{\hat{e}_{i}}\big|d\hat{s}
\,,\qquad i=1,\dots,4,\quad j=1,2\,,
\end{equation}
where $\mathbf{n}_{i}$ is the unit normal vector to the edge $e_{i}$, $\mathcal{F}$ is the mapping from the reference square in the natural coordinate system to the quadrilateral element $\mathscr{Q}$ in the physical space, $\mathcal{F}|_{\hat{e}_{i}}$ denotes its restriction to the edge $\hat{e}_{i}$, and $\big|\mathcal{F}|_{\hat{e}_{i}}\big|$ is the Jacobian of the transformation from the local coordinate system to the physical space. The same basis functions $\hat{h}_{j}$, $j=1,2$, defined in \eqref{eqn:2.7}, are used to compute the degrees of freedom in \eqref{eqn:2.15}. The coefficients of the shape functions are determined by enforcing the following conditions \citep{JAH25}:
\begin{equation}\label{eqn:2.16}
\varphi^{\mathscr{Q},e_{l}}_{m}\left(\mathbf{w}^{\mathscr{Q},e_{i}}_{j}\right)=
\begin{cases}
  1,&\text{if } l=i,~m=j\,, \\
  0,&\text{otherwise}\,.
\end{cases}
\end{equation}
The procedure described in Remark~\ref{rem:2.1} can also be used to efficiently compute the coefficients of these shape functions. To enforce the jump condition in \eqref{eqn:2.3}$_{2}$ for the normal component of $\mathbf{w}^{\mathscr{Q},e_{i}}_{j}$ across an edge $e$ shared by two quadrilaterals $\mathscr{Q}_{1}$ and $\mathscr{Q}_{2}$, we define the global shape functions $\mathbf{W}^{e}_{j}$, $j=1,2$, as follows:
\begin{equation}\label{eqn:2.17}
\mathbf{W}^{e}_{j}|_{\mathscr{Q}_{l}}=
\begin{cases}
  c^{l}_{i}\mathbf{w}^{\mathscr{Q}_{l},e_{i}}_{j}\,,&
  \text{if }\mathscr{Q}_{l}=\mathscr{Q}_{1},\mathscr{Q}_{2}\,,\\
  0,&\text{otherwise}\,,
\end{cases}
\end{equation}
where $e_{i}\in\partial\mathscr{Q}_{l}$ is the same edge with the global identifier $e$ and,
\begin{equation}\label{eqn:2.18}
c^{l}_{i}=\mathbf{n}^{e}\cdot\mathbf{n}^{l}_{i}\,.
\end{equation}
In \eqref{eqn:2.18}, $\mathbf{n}^{e}$ denotes the global unit normal vector associated with the edge $e$, while $\mathbf{n}^{l}_{i}$ denotes the local unit normal vector associated with the edge $e_{i}\in\partial\mathscr{Q}_{l}$. The global shape functions $\mathbf{W}^{e}_{j}$, $j=1,2$, associated with an edge $e$ shared by two quadrilaterals $\mathscr{Q}_{1}$ and $\mathscr{Q}_{2}$ are illustrated in Figure~\ref{fig:24}. As seen in the figure, the normal component of the corresponding vector field is single-valued across the common edge $e$.

We note that the degrees of freedom $\varphi^{\mathscr{Q},e_{i}}_{j}$, $i=1,\dots,4$, $j=1,2$, in \eqref{eqn:2.15} are sufficient to uniquely determine the shape functions used to interpolate the stress tensor. In particular, no additional degrees of freedom need to be defined on the quadrilateral $\mathscr{Q}$ itself. We also note that the order of the Lagrange shape functions used to interpolate the displacement field (bilinear in the present case) does not permit the use of higher-order shape functions for the stress tensor. If higher-order shape functions are used, the stiffness matrix relating the displacement degrees of freedom to the stress degrees of freedom fails to satisfy the conditions necessary to attain the minimum rank required for invertibility of the overall stiffness matrix. This point is discussed further in \S\ref{sec:3}.

\section{Compatible-Strain Mixed Finite Elements}
\label{sec:3}

In this section, we discuss two mixed finite element formulations applicable to the large-deformation analysis of compressible and incompressible two-dimensional solids. The first formulation is similar to those previously developed for simplicial elements \citep{JAH22,JAH25}. It is, however, modified to accommodate quadrilateral elements with different nodal arrangements and varying numbers of degrees of freedom. The second formulation is developed specifically for incompressible solids. In addition to the displacement, displacement gradient, and stress tensor used in the first formulation, the dilation and pressure are introduced as independent variables. Before presenting these formulations, it is necessary to discuss several preliminary aspects, including the assignment of degrees of freedom to nodes and pseudo-nodes, as well as the interpolation of the various fields.

\subsection{Degrees of freedom and field interpolations}
\label{sec:3.1}

\begin{figure}
\centering
\includegraphics[width=0.25\textwidth]{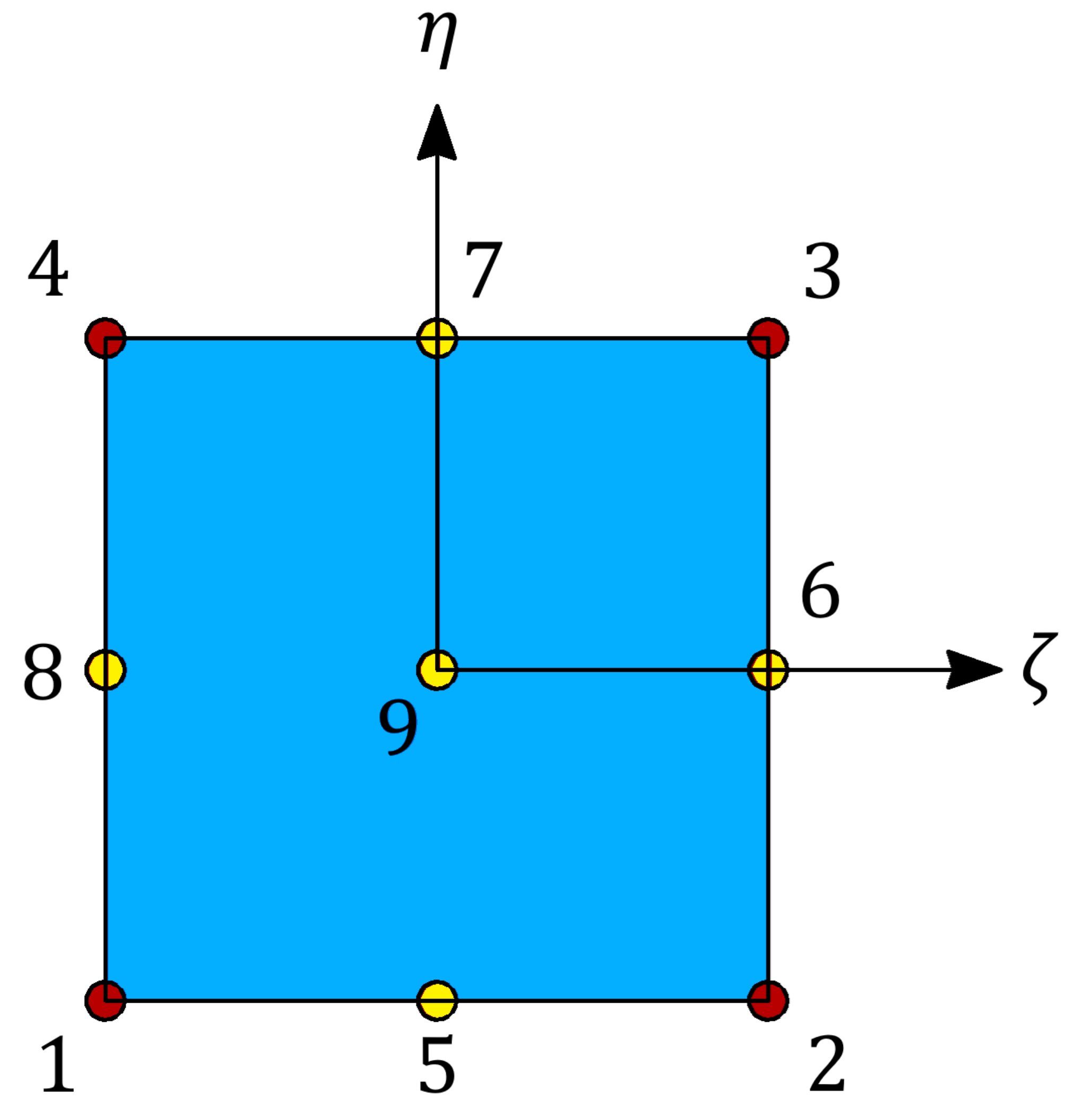}
\vskip 0.1in
\caption{Nodes, shown in red, and pseudo-nodes, shown in yellow, for the reference element in the natural coordinates $\zeta$ and $\eta$. The displacement degrees of freedom are defined at nodes $1$--$4$. The degrees of freedom associated with the displacement gradient are assigned to pseudo-nodes $5$--$9$. Pseudo-nodes $5$--$8$ also represent the degrees of freedom associated with the stress tensor.}
\label{fig:31}
\end{figure}

To formulate the element developed in this work, we consider four nodes and five pseudo-nodes, as shown in Figure~\ref{fig:31}. Pseudo-nodes $5$ to $8$ represent the degrees of freedom defined on the edges $e_{1}$ to $e_{4}$ of the quadrilateral $\mathscr{Q}$ (see Figure~\ref{fig:21}), while pseudo-node $9$ represents the degrees of freedom associated with $\mathscr{Q}$ itself. The displacement degrees of freedom are defined at nodes $1$ to $4$, resulting in a bilinear displacement field for the element. The degrees of freedom associated with the displacement gradient are assigned to pseudo-nodes $5$ to $9$. As discussed in \S\ref{sec:2}, it is not necessary to introduce additional degrees of freedom for the stress tensor on $\mathscr{Q}$ itself. Consequently, no stress degrees of freedom are assigned to pseudo-node $9$.

The displacement degrees of freedom $u_{i}$, $i=1,\dots,8$, are assigned to nodes $1$ to $4$. A total of $24$ degrees of freedom, denoted by $\alpha_{i}$, $i=1,\dots,24$, is used for the displacement gradient. Four degrees of freedom are assigned to each edge through its representative pseudo-node, and the remaining eight degrees of freedom are assigned to the quadrilateral $\mathscr{Q}$ through pseudo-node $9$. The stress tensor is represented by the degrees of freedom $\gamma_{i}$, $i=1,\dots,16$, with four degrees of freedom assigned to each edge. Therefore, nodes $1$ to $4$ carry two degrees of freedom each, while pseudo-nodes $5$ to $9$ carry eight degrees of freedom each. Table~\ref{tab:3.1} summarizes the association between the nodes and pseudo-nodes and the corresponding degrees of freedom assigned to them.
\begin{table}
\centering
\caption{Degrees of freedom for nodes $1$ to $4$ and pseudo-nodes $5$ to $9$ of the quadrilateral element $\mathscr{Q}$ shown in Figure~\ref{fig:31}.}
\label{tab:3.1}
\renewcommand{\arraystretch}{1.5}
\renewcommand{\tabcolsep}{0.2cm}
\begin{tabular}{c c c c}
\hline
Node/Pseudo-Node & \multicolumn{3}{c}{Degrees of freedom}                                         \\
\cline{2-4}
                 & Displacement  & Displacement Gradient            & Stress                      \\
\hline
1                & $u_{1},u_{2}$ &                                  &                             \\
2                & $u_{3},u_{4}$ &                                  &                             \\
3                & $u_{5},u_{6}$ &                                  &                             \\
4                & $u_{7},u_{8}$ &                                  &                             \\
5                &               & $\alpha_{1},\alpha_{2},\alpha_{3},\alpha_{4}$                  & $\gamma_{1},\gamma_{2},\gamma_{3},\gamma_{4}$     \\
6                &               & $\alpha_{5},\alpha_{6},\alpha_{7},\alpha_{8}$                  & $\gamma_{5},\gamma_{6},\gamma_{7},\gamma_{8}$     \\
7                &               & $\alpha_{9},\alpha_{10},\alpha_{11},\alpha_{12}$               &
$\gamma_{9},\gamma_{10},\gamma_{11},\gamma_{12}$  \\
8                &               & $\alpha_{13},\alpha_{14},\alpha_{15},\alpha_{16}$              &
$\gamma_{13},\gamma_{14},\gamma_{15},\gamma_{16}$ \\
9                &               & $\alpha_{17},\alpha_{18},\alpha_{19},\alpha_{20},\alpha_{21},\alpha_{22},\alpha_{23},\alpha_{24}$ &
                                                  \\
\hline
\end{tabular}
\end{table}

Similar to the first-order quadrilateral element, the displacement $\bm{\varphi}\left(\zeta,\eta\right):\left[-1,1\right]^{2}\rightarrow\mathbb{R}^{2}$ can be interpolated as follows \citep{BAT96,ZIE05,HUG00,BEL05,WRI08}:
\begin{equation}\label{eqn:3.1}
\bm{\varphi}\left(\zeta,\eta\right)=\sum_{i=1}^{4}N_{i}\left(\zeta,\eta\right)\mathbf{u}_{i}\,,
\end{equation}
where $N_{i}\left(\zeta,\eta\right):\left[-1,1\right]\rightarrow\mathbb{R},~i=1\dots,4$ are the Lagrange shape functions associated to nodes $1$ to $4$ and $\mathbf{u}^{\mathsf{T}}_{i}=\left[u_{2i-1}\quad u_{2i}\right]\in\mathbb{R}^{2},~i=1,\dots,4$ are the displacements of these nodes in $x$ and $y$ directions (see Figure~\ref{fig:21}). The displacement gradient $\mathbf{H}$, which belongs to the space $\mathcal{P}^{c}_{3}\left(\otimes^{2}\mathcal{Q}\right)$ of polynomial tensor fields of maximum degree $3$, is interpolated as:
\begin{equation}\label{eqn:3.2}
\mathbf{H}=\sum_{a=1}^{24}\alpha_{a}\mathbf{G}^{c}_{a}\,,
\end{equation}
where $\alpha_{a}$ are the degrees of freedom considered for displacement gradient (see Table~\ref{tab:3.1}). Based on the previous discussion regarding the displacement gradient degrees of freedom on edges $e_{i},~i=1,\dots,4$ and on the quadrilateral $\mathscr{Q}$, the index $a$ can be expressed as \citep{JAH25}:
\begin{equation}\label{eqn:3.3}
a=
\begin{cases}
4(i-1)+2(j-1)+l, & j\text{th shape function on edge } e_{i},\\
2(k-1)+l+16,     & k\text{th shape function on quadrilateral }\mathscr{Q}\,.
\end{cases}
\end{equation}
Here, $l=1,2$ refers to the interpolating tensor $\mathbf{G}^{c}_{a_{l}}$, which contains $\mathbf{V}^{e_{i}\mathsf{T}}_{j}$ or $\mathbf{V}^{\mathscr{Q}\,\mathsf{T}}_{j}$ (see \eqref{eqn:2.11} and \eqref{eqn:2.13}) on its $l$th row and the null matrix on the other row. More specifically, for the $j$th shape function on the edge $e_{i}$, we can compute $a_{1}$ and $a_{2}$ by substituting $l=1,2$ in \eqref{eqn:3.3}. The submatrices $\mathbf{G}^{c}_{a_{1}}$ and $\mathbf{G}^{c}_{a_{2}}$ can then be obtained as
\begin{equation}\label{eqn:3.4}
\mathbf{G}^{c}_{a_{1}}=
\begin{bmatrix}
  \mathbf{V}^{e_{i}\mathsf{T}}_{j} \\
  \mathbf{0}_{1\times2}
\end{bmatrix}\,,\qquad
\mathbf{G}^{c}_{a_{2}}=
\begin{bmatrix}
  \mathbf{0}_{1\times2} \\
  \mathbf{V}^{e_{i}\mathsf{T}}_{j}
\end{bmatrix}\,.
\end{equation}
The first Piola-Kirchhoff stress tensor $\mathbf{P}$ belongs to the space $\mathcal{P}^{d}_{2}\left(\otimes^{2}\mathcal{Q}\right)$ of polynomial tensor fields of maximum degree $2$ and it can be interpolated using the following relationship:
\begin{equation}\label{eqn:3.5}
\mathbf{P}=\sum_{b=1}^{16}\gamma_{b}\mathbf{G}^{d}_{b}\,,
\end{equation}
where $\gamma_{b}$ are the stress degrees of freedom on edges $e_{i},~i=1,\dots,4$. Noting that the stress degrees of freedom are defined on the edges only, the index $b$ can be expressed in the following form:
\begin{equation}\label{eqn:3.6}
b=4\left(i-1\right)+2\left(j-1\right)+l\,.
\end{equation}
Similarly, the parameter $l=1,2$ in \eqref{eqn:3.6} refers to the interpolating tensor $\mathbf{G}^{d}_{b_{l}}$, which contains $\mathbf{W}^{e_{i}\mathsf{T}}_{j}$ (see \eqref{eqn:2.17}) on its $l$th row and the null matrix on the other row. For the $j$th shape function on the edge $e_{i}$ we can compute $b_{1}$ and $b_{2}$ by substituting $l=1,2$ in \eqref{eqn:3.6} and obtain the submatrices $\mathbf{G}^{d}_{b_{1}}$ and $\mathbf{G}^{d}_{b_{2}}$ as follows:
\begin{equation}\label{eqn:3.7}
\mathbf{G}^{d}_{b_{1}}=
\begin{bmatrix}
  \mathbf{W}^{e_{i}\mathsf{T}}_{j} \\
  \mathbf{0}_{1\times2}
\end{bmatrix},\qquad
\mathbf{G}^{d}_{b_{2}}=
\begin{bmatrix}
  \mathbf{0}_{1\times2} \\
  \mathbf{W}^{e_{i}\mathsf{T}}_{j}
\end{bmatrix}\,.
\end{equation}

\subsection{Mixed formulation for compressible solids} \label{sec:3.2}

In the following, we consider an elastic body occupying a region $\mathcal{B}\subset\mathbb{R}^2$ in the reference configuration, with $\partial_{\tau}\mathcal{B}\subset\partial\mathcal{B}$ denoting the portion of the boundary on which tractions are prescribed. Material points in the reference configuration are identified by the coordinates $X^{i}$, $i=1,2$, while their positions in the current configuration are described by the coordinates $x^{i}$, $i=1,2$. The displacement vector $\bm{\varphi}=\mathbf{x}-\mathbf{X}$, the displacement gradient $\mathbf{H}$, and the first Piola--Kirchhoff stress tensor $\mathbf{P}$ are taken as independent variables in the Hu--Washizu functional, which is given by \citep{JAH25}\footnote{The Hu--Washizu functional has its roots in earlier mixed variational principles in elasticity. \citet{Hellinger1914} formulated one of the earliest mixed principles by treating stresses as independent variables, \citet{Reissner1950} later developed a two-field variational principle involving independent stress and displacement fields. The three-field formulation, in which displacement, strain, and stress are treated as independent fields, was introduced independently by \citet{Hu1955} and \citet{Washizu1955}. Washizu subsequently systematized and popularized this approach in his monograph \citep{Washizu1982}, and the resulting Hu--Washizu principle has since become one of the fundamental variational principles underlying mixed finite element formulations.}
\begin{equation}\label{eqn:3.8}
\Pi_{1}\left(\bm{\varphi},\mathbf{H},\mathbf{P}\right)=
\int_{\mathcal{B}}\widehat{W}\left(\mathbf{C}\right)dV+
\int_{\mathcal{B}}\bm{\tau}:\left(\nabla\bm{\varphi}-\mathbf{h}\right)dV-
\int_{\mathcal{B}}\rho_{0}\left(\mathbf{b}\cdot\bm{\varphi}\right)dV-
\int_{\partial_{\tau}\mathcal{B}}\mathbf{t}\cdot\bm{\varphi}\,dA\,,
\end{equation}
where $\widehat{W}$ is the stored energy function, $\mathbf{C}=\mathbf{F}^{\mathsf{T}}\mathbf{F}$ is the right Cauchy--Green tensor, $\mathbf{F}=\mathbf{I}+\mathbf{H}$ is the deformation gradient, $\mathbf{h}=\mathbf{H}\mathbf{F}^{-1}$ is the push-forward of the displacement gradient to the current configuration, and $\bm{\tau}=\mathbf{P}\mathbf{F}^{\mathsf{T}}$ is the Kirchhoff stress tensor. Furthermore, $\rho_{0}$ denotes the material mass density, $\mathbf{b}$ the body force, and $\mathbf{t}$ the traction prescribed on the portion of the boundary where tractions are specified. The differential operators $\nabla_{0}=\frac{\partial}{\partial\mathbf{X}}$ and $\nabla=\left(\frac{\partial}{\partial\mathbf{X}}\right)\mathbf{F}^{-1}$ denote differentiation with respect to the material and spatial configurations, respectively. We note that the Kirchhoff stress tensor $\bm{\tau}$ acts as a Lagrange multiplier enforcing the constraint $\nabla\bm{\varphi}=\mathbf{h}$, which is equivalent to $\nabla_{0}\bm{\varphi}=\mathbf{H}$.

The governing equations for the body $\mathcal{B}$ under the applied loads are obtained by setting the first variation of the functional in \eqref{eqn:3.8} to zero. We consider the independent variations $\delta\bm{\varphi}$, $\delta\mathbf{H}$, and $\delta\mathbf{P}$ of $\bm{\varphi}$, $\mathbf{H}$, and $\mathbf{P}$, respectively. Denoting the derivative of the functional $\Pi_{1}$ in the direction of the variation $\delta\bm{\varphi}$ by $D\Pi_{1}\cdot\delta\bm{\varphi}$, the variations of $\Pi_{1}$ corresponding to the aforementioned independent variations are given by $(\delta_{\delta\bm{\varphi}}\Pi_{1},\delta_{\delta\mathbf{H}}\Pi_{1},\delta_{\delta\mathbf{P}}\Pi_{1})=(D\Pi_{1}\cdot\delta\bm{\varphi},D\Pi_{1}\cdot\delta\mathbf{H},D\Pi_{1}\cdot\delta\mathbf{P})$. These three variations are calculated as follows:
\begin{equation}\label{eqn:3.9}
\begin{aligned}
D\Pi_{1}\cdot\delta\bm{\varphi} &=
\int_{\mathcal{B}}\bm{\tau}:\nabla\delta\bm{\varphi}\,dV-
\int_{\mathcal{B}}\rho_{0}\left(\mathbf{b}\cdot\delta\bm{\varphi}\right)dV-
\int_{\partial_{\tau}\mathcal{B}}\mathbf{t}\cdot\delta\bm{\varphi}\,dA=0\,, \\
D\Pi_{1}\cdot\delta\mathbf{H} &=
\int_{\mathcal{B}}\bm{\widehat{\tau}}:\delta\mathbf{h}^{s}\,dV-
\int_{\mathcal{B}}\bm{\tau}:\delta\mathbf{h}\,dV=0\,, \\
D\Pi_{1}\cdot\delta\mathbf{P} &=
\int_{\mathcal{B}}\delta\bm{\tau}:\nabla\bm{\varphi}\,dV-
\int_{\mathcal{B}}\delta\bm{\tau}:\mathbf{h}\,dV=0\,,
\end{aligned}
\end{equation}
where $\bm{\widehat{\tau}}=2\mathbf{F}\frac{\partial\widehat{W}}{\partial\mathbf{C}}\mathbf{F}^{\mathsf{T}}$ is the Kirchhoff stress tensor associated with the stored energy function $\widehat{W}$, and $\delta\mathbf{h}^s=\frac{1}{2}\left(\delta\mathbf{h}+\delta\mathbf{h}^{\mathsf{T}}\right)$. We note that Eqs.~\eqref{eqn:3.9} are nonlinear and, therefore, must be solved iteratively. Following \citet{JAH25}, we introduce $\Delta\bm{\varphi}$, $\Delta\mathbf{H}$, and $\Delta\mathbf{P}$ as another set of independent variations of $\bm{\varphi}$, $\mathbf{H}$, and $\mathbf{P}$, respectively, and linearize Eqs.~\eqref{eqn:3.9} as follows:
\begin{equation}\label{eqn:3.10}
\begin{aligned}
&\int_{\mathcal{B}}\Delta\bm{\tau}:\nabla\delta\bm{\varphi}\,dV=
\int_{\mathcal{B}}\rho_{0}\left(\mathbf{b}\cdot\delta\bm{\varphi}\right)dV+
\int_{\partial_{\tau}\mathcal{B}}\mathbf{t}\cdot\delta\bm{\varphi}\,dA-
\int_{\mathcal{B}}\bm{\tau}:\nabla\delta\bm{\varphi}\,dV\,, \\
&\int_{\mathcal{B}}\delta\mathbf{h}^{s}:\left(\mathbb{C}:\Delta\mathbf{h}^{s}\right)dV+
\int_{\mathcal{B}}\left(\Delta\mathbf{h}\,\bm{\widehat{\tau}}\right):\delta\mathbf{h}\,dV-
\int_{\mathcal{B}}\Delta\bm{\tau}:\delta\mathbf{h}\,dV=
-\left(\int_{\mathcal{B}}\bm{\widehat{\tau}}:\delta\mathbf{h}^{s}\,dV
-\int_{\mathcal{B}}\bm{\tau}:\delta\mathbf{h}\,dV\right)\,, \\
&\int_{\mathcal{B}}\delta\bm{\tau}:\nabla\Delta\bm{\varphi}\,dV-
\int_{\mathcal{B}}\delta\bm{\tau}:\Delta\mathbf{h}\,dV=
-\left(\int_{\mathcal{B}}\delta\bm{\tau}:\nabla\bm{\varphi}\,dV-
\int_{\mathcal{B}}\delta\bm{\tau}:\mathbf{h}\,dV\right)\,.
\end{aligned}
\end{equation}
In these equations, $\Delta\bm{\tau}=\Delta\mathbf{P}\,\mathbf{F}^{\mathsf{T}}$, $\boldsymbol{\mathbb{C}}$ is the spatial fourth-order elasticity tensor, $\Delta\mathbf{h}=\Delta\mathbf{H}\,\mathbf{F}^{-1}$, and $\Delta\mathbf{h}^{s}$ denotes the symmetric part of $\Delta\mathbf{h}$.

From Eqs.~\eqref{eqn:3.2} and \eqref{eqn:3.5}, it follows that the tensor $\mathbf{h}$ and the Kirchhoff stress tensor $\bm{\tau}$ can be interpolated as:
\begin{equation}\label{eqn:3.11}
\mathbf{h}=\sum_{a=1}^{24}\alpha_{a}\mathbf{g}^{c}_{a}\,,\qquad
\bm{\tau}=\sum_{b=1}^{16}\gamma_{b}\mathbf{g}^{d}_{b}\,,
\end{equation}
where $\mathbf{g}^{c}_{a}=\mathbf{G}^{c}_{a}\mathbf{F}^{-1}$ and $\mathbf{g}^{d}_{b}=\mathbf{G}^{d}_{b}\mathbf{F}^{\mathsf{T}}$. In order to express \eqref{eqn:3.10} in matrix form, we represent the components of the tensors $\mathbf{h}$ and $\bm{\tau}$ in vector form as:
\begin{equation}\label{eqn:3.12}
\bm{\varepsilon}=
\begin{bmatrix}
  h_{11} \\
  h_{12} \\
  h_{21} \\
  h_{22}
\end{bmatrix}\,,\qquad
\bm{\Gamma}=
\begin{bmatrix}
  \tau_{11} \\
  \tau_{12} \\
  \tau_{21} \\
  \tau_{22}
\end{bmatrix}\,.
\end{equation}
In order to be consistent with the definition of the vectors $\bm{\varepsilon}$ and $\bm{\Gamma}$, we define the matrices $\mathbf{\bar{g}}^{c}$ and $\mathbf{\bar{g}}^{d}$ as follows:
\begin{equation}\label{eqn:3.13}
\mathbf{\bar{g}}^{c}=
\begin{bmatrix}
  \left(\mathbf{g}^{c}_{1}\right)_{11}&\left(\mathbf{g}^{c}_{2}\right)_{11}&\cdots&\left(\mathbf{g}^{c}_{24}\right)_{11}\\
  \left(\mathbf{g}^{c}_{1}\right)_{12}&\left(\mathbf{g}^{c}_{2}\right)_{12}&\cdots&\left(\mathbf{g}^{c}_{24}\right)_{12}\\
  \left(\mathbf{g}^{c}_{1}\right)_{21}&\left(\mathbf{g}^{c}_{2}\right)_{21}&\cdots&\left(\mathbf{g}^{c}_{24}\right)_{21}\\
  \left(\mathbf{g}^{c}_{1}\right)_{22}&\left(\mathbf{g}^{c}_{2}\right)_{22}&\cdots&\left(\mathbf{g}^{c}_{24}\right)_{22}
\end{bmatrix}_{4\times24}\,,\qquad
\mathbf{\bar{g}}^{d}=
\begin{bmatrix}
  \left(\mathbf{g}^{d}_{1}\right)_{11}&\left(\mathbf{g}^{d}_{2}\right)_{11}&\cdots&\left(\mathbf{g}^{d}_{16}\right)_{11}\\
  \left(\mathbf{g}^{d}_{1}\right)_{12}&\left(\mathbf{g}^{d}_{2}\right)_{12}&\cdots&\left(\mathbf{g}^{d}_{16}\right)_{12}\\
  \left(\mathbf{g}^{d}_{1}\right)_{21}&\left(\mathbf{g}^{d}_{2}\right)_{21}&\cdots&\left(\mathbf{g}^{d}_{16}\right)_{21}\\
  \left(\mathbf{g}^{d}_{1}\right)_{22}&\left(\mathbf{g}^{d}_{2}\right)_{22}&\cdots&\left(\mathbf{g}^{d}_{16}\right)_{22}
\end{bmatrix}_{4\times16}\,.
\end{equation}
The $a$th column of the matrix $\bar{\mathbf{g}}^{c}$ contains the components of the matrix $\mathbf{g}^{c}_{a}$ defined in \eqref{eqn:3.11}$_{1}$. Similarly, the $b$th column of the matrix $\bar{\mathbf{g}}^{d}$ contains the components of the matrix $\mathbf{g}^{d}_{b}$ defined in \eqref{eqn:3.11}$_{2}$. Using \eqref{eqn:3.12}--\eqref{eqn:3.13}, \eqref{eqn:3.11} can be written in the following matrix form:
\begin{equation}\label{eqn:3.15}
\bm{\varepsilon}=\mathbf{\bar{g}}^{c}\bm{\alpha}\,,\qquad
\bm{\Gamma}=\mathbf{\bar{g}}^{d}\bm{\gamma}\,,
\end{equation}
where the vectors $\bm{\alpha}$ and $\bm{\gamma}$ list the degrees of freedom $\alpha_{a},~a=1,\ldots,24$ and $\gamma_{b},~b=1,\ldots,16$, respectively. In accordance with the definition of the vector $\bm{\varepsilon}$ in \eqref{eqn:3.12}$_{1}$, we define the gradient of the displacement field $\bm{\varphi}$ as follows:
\begin{equation}\label{eqn:3.16}
\bm{\epsilon}^{\mathsf{T}}=
\begin{bmatrix}
  \left(\nabla\bm{\varphi}\right)_{11} & \left(\nabla\bm{\varphi}\right)_{12} &
  \left(\nabla\bm{\varphi}\right)_{21} & \left(\nabla\bm{\varphi}\right)_{22}
\end{bmatrix}\,.
\end{equation}
Noting that the displacement field is interpolated using \eqref{eqn:3.1}, one can define the matrices $\mathbf{N}$ and $\mathbf{B}$ as:
\begin{equation}\label{eqn:3.17}
\mathbf{N}=
\begin{bmatrix}
  N_{1} & 0     & \ldots & N_{4} & 0     \\
  0     & N_{1} & \ldots & 0     & N_{4}
\end{bmatrix}_{2\times8}\,,\qquad
\mathbf{B}=
\begin{bmatrix}
  \nabla N_{1}          & \mathbf{0}_{2\times1} & \ldots                &
  \nabla N_{4}          & \mathbf{0}_{2\times1} \\
  \mathbf{0}_{2\times1} & \nabla N_{1}          & \ldots                &
  \mathbf{0}_{2\times1} & \nabla N_{4}
\end{bmatrix}_{4\times8}\,,
\end{equation}
where the gradients of shape functions in the spatial configuration are obtained from $\nabla N_{i}=\left(\nabla_{0}N_{i}\right)\mathbf{F}^{-1}$. Using \eqref{eqn:3.17} we can interpolate the displacement field $\bm{\varphi}$ and its gradient $\bm{\epsilon}$ via the following equations:
\begin{equation}\label{eqn:3.19}
\bm{\varphi}=\mathbf{N}\mathbf{u}\,,\qquad
\bm{\epsilon}=\mathbf{B}\mathbf{u}\,,
\end{equation}
where the vector $\mathbf{u}$ lists the displacement degrees of freedom $u_{i},~i=1,\dots,8$.

Using Eqs.~\eqref{eqn:3.15}, \eqref{eqn:3.16}, and \eqref{eqn:3.19}, we represent $\delta\bm{\varphi}$, $\nabla\delta\bm{\varphi}$, $\delta\mathbf{h}$, and $\delta\bm{\tau}$ in the following form
\begin{equation}\label{eqn:3.20}
\delta\bm{\varphi}=\mathbf{N}\,\delta\mathbf{u}\,,\qquad
\delta\bm{\epsilon}=\mathbf{B}\,\delta\mathbf{u}\,,\qquad
\delta\bm{\varepsilon}=\mathbf{\bar{g}}^{c}\delta\bm{\alpha},\qquad
\delta\bm{\Gamma}=\mathbf{\bar{g}}^{d}\delta\bm{\gamma}\,.
\end{equation}
Similarly, $\Delta\bm{\varphi}$, $\nabla\Delta\bm{\varphi}$, $\Delta\mathbf{h}$ and $\Delta\bm{\tau}$ are expressed in the following form
\begin{equation}\label{eqn:3.21}
\Delta\bm{\varphi}=\mathbf{N}\,\Delta\mathbf{u}\,,\qquad
\Delta\bm{\epsilon}=\mathbf{B}\,\Delta\mathbf{u},\qquad
\Delta\bm{\varepsilon}=\mathbf{\bar{g}}^{c}\Delta\bm{\alpha}\,,\qquad
\Delta\bm{\Gamma}=\mathbf{\bar{g}}^{d}\Delta\bm{\gamma}\,.
\end{equation}
Noting that the vectors $\delta\mathbf{u}$, $\delta\bm{\alpha}$ and $\delta\bm{\gamma}$ are arbitrary, it can be shown that \eqref{eqn:3.10} reduce to the following matrix equations \citep{JAH25}:
\begin{align}\label{eqn:3.22}
\mathbf{K}^{u\gamma}_{\iota}\Delta\bm{\gamma}=
\mathbf{F}_{b}+\mathbf{F}_{t}-\mathbf{F}^{u}_{i}\,,\qquad
\mathbf{K}^{\alpha\alpha}_{m}\Delta\bm{\alpha}-
\mathbf{K}^{\alpha\gamma}_{\iota}\Delta\bm{\gamma}=-\mathbf{F}^{\alpha}_{i}\,,\qquad
\mathbf{K}^{\gamma u}_{\iota}\Delta\mathbf{u}-
\mathbf{K}^{\gamma\alpha}_{\iota}\Delta\bm{\alpha}=-\mathbf{F}^{\gamma}_{i}\,,
\end{align}
where the stiffness matrices are defined as:
\begin{equation} \label{eqn:3.23}
\begin{aligned}
& \mathbf{K}^{u\gamma}_{\iota}=\mathbf{K}^{\gamma u\,\mathsf{T}}_{\iota}=
\int_{\mathcal{B}}\mathbf{B}^{\mathsf{T}}\mathbf{\bar{g}}^{d}dV,\qquad
\mathbf{K}^{\alpha\alpha}_{m}=
\int_{\mathcal{B}}\mathbf{\bar{g}}^{c\,\mathsf{T}}\mathbf{D}\,\mathbf{\bar{g}}^{c}dV+
\int_{\mathcal{B}}\mathbf{\bar{g}}^{c\,\mathsf{T}}\mathbf{\widehat{T}}\,\mathbf{\bar{g}}^{c}dV\,,\\
&\mathbf{K}^{\alpha\gamma}_{\iota}=\mathbf{K}^{\gamma\alpha\,\mathsf{T}}_{\iota}=
\int_{\mathcal{B}}\mathbf{\bar{g}}^{c\,\mathsf{T}}\mathbf{\bar{g}}^{d}dV\,.
\end{aligned}
\end{equation}
The matrix $\mathbf{D}$ in \eqref{eqn:3.23}$_{2}$ represents the elasticity tensor $\mathbb{C}$ in matrix form and in view of the definition of the vector $\bm{\varepsilon}$ in \eqref{eqn:3.12} it has the following representation:
\begin{equation}\label{eqn:3.24}
\mathbf{D}=
\begin{bmatrix}
\mathbb{C}_{1111}&\mathbb{C}_{1112}&\mathbb{C}_{1121}&\mathbb{C}_{1122} \\
\mathbb{C}_{1211}&\mathbb{C}_{1212}&\mathbb{C}_{1221}&\mathbb{C}_{1222} \\
\mathbb{C}_{2111}&\mathbb{C}_{2112}&\mathbb{C}_{2121}&\mathbb{C}_{2122} \\
\mathbb{C}_{2211}&\mathbb{C}_{2212}&\mathbb{C}_{2221}&\mathbb{C}_{2222}
\end{bmatrix}\,.
\end{equation}
Furthermore, the matrix $\mathbf{\widehat{T}}$ in \eqref{eqn:3.23}$_{2}$ has the following representation
\begin{equation}\label{eqn:3.25}
\mathbf{\widehat{T}}=
\begin{bmatrix}
  \bm{\widehat{\tau}}   & \mathbf{0}_{4\times4} \\
  \mathbf{0}_{4\times4} & \bm{\widehat{\tau}}
\end{bmatrix}\,.
\end{equation}
The vectors $\mathbf{F}_{b}$ and $\mathbf{F}_{t}$ in \eqref{eqn:3.22} are the body and traction force vectors, which are defined, respectively, as
\begin{equation}\label{eqn:3.26}
\mathbf{F}_{b}=\int_{\mathcal{B}}\rho_{0}\left(\mathbf{N}^{\mathsf{T}}\mathbf{b}\right)dV,\qquad
\mathbf{F}_{t}=\sum_{b=1}^{n_{b}}\int_{\partial_{\tau}M_{b}}\mathbf{N}^{\mathsf{T}}_{b}\mathbf{t}\,dA\,,
\end{equation}
where the matrix $\mathbf{N}_{b}$ is defined analogously to the matrix $\mathbf{N}$ in \eqref{eqn:3.17} and is used to interpolate the displacement field on the portion of the boundary where tractions are prescribed. From \eqref{eqn:3.10}, we observe that the internal load vectors $\mathbf{F}^{u}_{i}$, $\mathbf{F}^{\alpha}_{i}$, and $\mathbf{F}^{\gamma}_{i}$ can be defined as follows:
\begin{equation}\label{eqn:3.27}
\mathbf{F}^{u}_{i}=
\int_{\mathcal{B}}\mathbf{B}^{\mathsf{T}}\bm{\Gamma}\,dV,\qquad
\mathbf{F}^{\alpha}_{i}=
\int_{\mathcal{B}}\mathbf{\bar{g}}^{c\,\mathsf{T}}\left(\bm{\widehat{\Gamma}}-\bm{\Gamma}\right)dV,\qquad
\mathbf{F}^{\gamma}_{i}=
\int_{\mathcal{B}}\mathbf{\bar{g}}^{d\,\mathsf{T}}\left(\bm{\epsilon}-\bm{\varepsilon}\right)dV\,.
\end{equation}
Eqs.~\eqref{eqn:3.22} can be conveniently cast into the following form
\begin{equation}\label{eqn:3.28}
\mathbf{K}^{\mathscr{Q}}_{t}\Delta\mathbf{U}^{\mathscr{Q}}=\mathbf{F}^{\mathscr{Q}}_{e}-\mathbf{F}^{\mathscr{Q}}_{i}\,,
\end{equation}
where the superscript $\mathscr{Q}$ implies that the equation corresponds to a given quadrilateral $\mathscr{Q}\in\mathcal{Q}$. The tangent stiffness matrix in \eqref{eqn:3.28} is expressed as
\begin{equation}\label{eqn:3.29}
\mathbf{K}^{\mathscr{Q}}_{t}=
\begin{bmatrix}
  \mathbf{0}_{8 \times8}        &  \mathbf{0}_{8\times24}           &  \mathbf{K}^{u\gamma}_{\iota}      \\
  \mathbf{0}_{24\times8}        &  \mathbf{K}^{\alpha\alpha}_{m}    & -\mathbf{K}^{\alpha\gamma}_{\iota} \\
  \mathbf{K}^{\gamma u}_{\iota} & -\mathbf{K}^{\gamma\alpha}_{\iota}&  \mathbf{0}_{16\times16}
\end{bmatrix}\,.
\end{equation}
On the other hand, the vectors $\Delta\mathbf{U}^\mathscr{Q}$, $\mathbf{F}^\mathscr{Q}_{e}$ and $\mathbf{F}^\mathscr{Q}_{i}$ have the following representations:
\begin{equation}\label{eqn:3.30}
\Delta\mathbf{U}^{\mathscr{Q}}=
\begin{bmatrix}
  \Delta\mathbf{u}  \\
  \Delta\bm{\alpha} \\
  \Delta\bm{\gamma}
\end{bmatrix}\,,\qquad
\mathbf{F}^{\mathscr{Q}}_{e}=
\begin{bmatrix}
  \mathbf{F}_{b}+\mathbf{F}_{t} \\
  \mathbf{0}_{24\times1}        \\
  \mathbf{0}_{16\times1}
\end{bmatrix}\,,\qquad
\mathbf{F}^{\mathscr{Q}}_{i}=
\begin{bmatrix}
  \mathbf{F}^{u}_{i}      \\
  \mathbf{F}^{\alpha}_{i} \\
  \mathbf{F}^{\gamma}_{i}
\end{bmatrix}\,.
\end{equation}
The global stiffness matrix, and the global external and internal load vectors can be obtained by assembling the stiffness matrix, external and internal load vectors for each quadrilateral $\mathscr{Q}$ as specified by the following relations:
\begin{equation}\label{eqn:3.31}
\mathbb{K}_{t}=\underset{\mathscr{Q}\in\mathcal{Q}}{\mathbf{\scalebox{1.5}{A}}}\mathbf{K}^{\mathscr{Q}}_{t}\,,\qquad
\mathbb{F}_{e}=\underset{\mathscr{Q}\in\mathcal{Q}}{\mathbf{\scalebox{1.5}{A}}}\mathbf{F}^{\mathscr{Q}}_{e}\,,\qquad
\mathbb{F}_{i}=\underset{\mathscr{Q}\in\mathcal{Q}}{\mathbf{\scalebox{1.5}{A}}}\mathbf{F}^{\mathscr{Q}}_{i}\,,
\end{equation}
where $\mathbf{\scalebox{1.5}{A}}$ is the assembly operator. The vector $\Delta\mathbb{U}$, obtained by solving the assembled system of equations, is used to update the global degrees of freedom via the relation $\mathbb{U}_{i+1} = \mathbb{U}_{i} + \Delta\mathbb{U}_{i}$. The subscript $i$ in this relation refers to the iteration number. The iterative procedure continues until the norm of $\Delta\mathbb{U}_{i}$ is negligible within a given tolerance value.

\subsection{Mixed formulation for incompressible solids} \label{sec:3.3}

For incompressible solids, we augment the independent variables $(\bm{\varphi},\mathbf{H},\mathbf{P})$ of the Hu--Washizu functional by introducing the dilation $\theta$ and the pressure $p$. We assume that the stored energy function $\widehat{W}$ is additively decomposed as $\widehat{W}=\widehat{W}(\bar{\mathbf{C}})+U(\theta)$, where $\bar{\mathbf{C}}=J^{-\frac{2}{3}}\mathbf{C}$. Under these assumptions, the Hu--Washizu functional takes the following form:
\begin{equation}\label{eqn:3.32}
\begin{aligned}
\Pi_{2}\left(\bm{\varphi},\mathbf{H},\mathbf{P},\theta,p\right)&=
\int_{\mathcal{B}}\left[\widehat{W}\left(\mathbf{\bar{C}}\right)+U\left(\theta\right)\right]dV+
\int_{\mathcal{B}}\bm{\tau}:\left(\nabla\bm{\varphi}-\mathbf{h}\right)dV+
\int_{\mathcal{B}}p\left(J-\theta\right)dV \\
& \quad -\int_{\mathcal{B}}\rho_{0}\left(\mathbf{b}\cdot\bm{\varphi}\right)dV-
\int_{\partial_{\tau}\mathcal{B}}\mathbf{t}\cdot\bm{\varphi}\,dA\,.
\end{aligned}
\end{equation}
In \eqref{eqn:3.32}, $J$ is the Jacobian of deformation.\footnote{Suppose the local coordinates $\{X^A\}$ and $\{x^a\}$ are used for the reference configuration $\mathcal{B}$ and the Euclidean ambient space $\mathcal{S}$. The metric of the ambient space is denoted as $\mathbf{g}$, which induces the flat metric $\mathbf{G}=\mathbf{g}\big|_{\mathcal{B}}$ on the reference configuration. The Jacobian of deformation has the following expression: $J=\sqrt{\frac{\det\mathbf{g}}{\det\mathbf{G}}}\,\det\mathbf{F}$ \citep{MarsdenHughes1983}. When Cartesian coordinates are used for both the reference and current configurations, this expression is simplified to $J=\det\mathbf{F}$.} Considering the independent variations $\delta\bm{\varphi}$, $\delta\mathbf{H}$, $\delta\mathbf{P}$, $\delta\theta$ and $\delta p$, respectively, in $\bm{\varphi}$, $\mathbf{H}$, $\mathbf{P}$, $\theta$ and $p$, and setting to zero the corresponding variations of the functional $\Pi_{2}$ leads to the set of governing equations. It can be verified that $\delta_{\delta\bm{\varphi}}\Pi_{2}=\delta_{\delta\bm{\varphi}}\Pi_{1}$ and $\delta_{\delta\mathbf{P}}\Pi_{2}=\delta_{\delta\mathbf{P}}\Pi_{1}$, and hence Eqs.~\eqref{eqn:3.9}$_{1}$ and \eqref{eqn:3.9}$_{3}$ remain valid for the case of incompressible solids. However, using the variation $\delta_{\delta\mathbf{H}}\Pi_{2}$ the Kirchhoff stress tensor $\widehat{\bm{\tau}}$ in \eqref{eqn:3.9}$_{2}$ takes the following form:
\begin{equation}\label{eqn:3.33}
\widehat{\bm{\tau}}=J^{-\frac{2}{3}}\left(\operatorname{dev}\:\bm{\bar{\tau}}\right)+J\,p\,\mathbf{I}\,,
\end{equation}
where $\bm{\bar{\tau}}=2\mathbf{F}\frac{\partial\widehat{W}}{\partial\mathbf{\bar{C}}}\mathbf{F}^{\mathsf{T}}$ and $\mathbf{I}$ is the spatial identity tensor of second-order. The deviator operator $\operatorname{dev}$ when applied to a second-order tensor $\mathbf{a}$ gives $\operatorname{dev}\:\mathbf{a}=\mathbf{a}-\frac{1}{3}\left(\mathrm{tr}\:\mathbf{a}\right)\mathbf{I}$, where the factor $\frac{1}{3}$, rather than $\frac{1}{2}$, is used because the plane-strain formulation is derived from a three-dimensional constitutive model and the stress tensor remains three-dimensional with a generally nonzero out-of-plane component. Computing the variations $\delta_{\delta\theta}\Pi_{2}$ and $\delta_{\delta p}\Pi_{2}$ and setting them to zero leads to the following equations:
\begin{equation}\label{eqn:3.34}
D\Pi_{2}\cdot\delta\theta=\int_{\mathcal{B}}\left[U'\left(\theta\right)-p\right]\delta\theta\;dV=0\,,\qquad
D\Pi_{2}\cdot\delta P=\int_{\mathcal{B}}\left(J-\theta\right)\delta p\;dV=0\,.
\end{equation}

The set of equations \eqref{eqn:3.9} can be linearized following the approach of \citet{JAH25}. It is important to note that the pressure $p$ is condensed out at the element level using \eqref{eqn:3.34}. We define the residuals $R_{1}\left(\delta\bm{\varphi}\right)=\delta_{\delta\bm{\varphi}}\Pi_{2}$, $R_{2}\left(\delta\mathbf{H}\right)=\delta_{\delta\mathbf{H}}\Pi_{2}$, and $R_{3}\left(\delta\mathbf{P}\right)=\delta_{\delta\mathbf{P}}\Pi_{2}$. Then, considering $\Delta\bm{\varphi}$, $\Delta\mathbf{H}$, and $\Delta\mathbf{P}$ as another set of independent variations of $\bm{\varphi}$, $\mathbf{H}$, and $\mathbf{P}$, respectively, we linearize \eqref{eqn:3.9} as follows:
\begin{equation}\label{eqn:3.35}
R_{1}\left(\delta\bm{\varphi}\right)+\delta R_{1}\left(\delta\bm{\varphi}\right)=0\,,\qquad
R_{2}\left(\delta\mathbf{H}\right)+\delta R_{2}\left(\delta\mathbf{H}\right)=0\,,\qquad
R_{3}\left(\delta\mathbf{P}\right)+\delta R_{3}\left(\delta\mathbf{P}\right)=0\,.
\end{equation}
It is straightforward to show that,
\begin{equation} \label{eqn:3.36}
\begin{aligned}
\delta R_{1}\left(\delta\bm{\varphi}\right)&=
DR_{1}\left(\delta\bm{\varphi}\right)\cdot\Delta\mathbf{P}
=\int_{\mathcal{B}}\Delta\bm{\tau}:\nabla\delta\bm{\varphi}\,dV, \\
\delta R_{2}\left(\delta\mathbf{H}\right)&=
DR_{2}\left(\delta\mathbf{H}\right)\cdot\Delta\mathbf{H}+DR_{2}\left(\delta\mathbf{H}\right)\cdot\Delta\mathbf{P} \\
&=\int_{\mathcal{B}}\delta\mathbf{h}^{s}:\left[\boldsymbol{\mathbb{C}}^{d}+pJ\left(\mathbf{I}\otimes\mathbf{I}-2\mathcal{I}\right)\right]:\Delta\mathbf{h}^{s}dV+
\int_{\mathcal{B}}\left(\Delta\mathbf{h}\,\bm{\widehat{\tau}}\right):\delta\mathbf{h}\,dV \\
&\quad +\int_{\mathcal{B}}\left(Dp\cdot\Delta\mathbf{H}\right)J\left(\mathrm{tr}\:\delta\mathbf{h}^{s}\right)\,dV-
\int_{\mathcal{B}}\Delta\bm{\tau}:\delta\mathbf{h}\,dV, \\
\delta R_{3}\left(\delta\mathbf{P}\right)&=
DR_{3}\left(\delta\mathbf{P}\right)\cdot\Delta\bm{\varphi}+DR_{3}\left(\delta\mathbf{P}\right)\cdot\Delta\mathbf{H}
=\int_{\mathcal{B}}\delta\bm{\tau}:\nabla\Delta\bm{\varphi}\,dV-
\int_{\mathcal{B}}\delta\bm{\tau}:\Delta\mathbf{h}\,dV\,,
\end{aligned}
\end{equation}
where $\mathbb{C}^{d}=\phi_{\star}\left[\frac{\partial^{2}\widehat{W}}{\partial\mathbf{C}\partial\mathbf{C}}\right]$\footnote{For the material fourth-order tensor $\mathbf{\mathcal{C}}$, its push-forward $\mathbb{C}=\phi_{\star}\left[\mathbf{\mathcal{C}}\right]$ to the current configuration has components $\mathbb{C}^{abcd}=F^{a}{}_{A}F^{b}{}_{B}F^{c}{}_{C}F^{d}{}_{D}\mathcal{C}^{ABCD}$, where $F^{a}{}_{A}$, $a,A=1,2$, are the components of the deformation gradient $\mathbf{F}$, and summation over repeated indices is assumed.} and $\mathcal{I}$ is the fourth-order identity tensor.\footnote{If the metric of the ambient space $\mathcal{S}$ is denoted by $\mathbf{g}$, then the components of the second- and fourth-order identity tensors with respect to a local chart $\{x^a\}$ are $I^{ab}=g^{ab}$ and $\mathcal{I}^{abcd}=\frac{1}{2}\left(g^{ac}g^{bd}+g^{ad}g^{bc}\right)$, respectively.} In the third integral of \eqref{eqn:3.36}$_{2}$, the pressure $p$ is treated as an independent variable and therefore does not admit a variation with respect to $\Delta\mathbf{H}$. However, as discussed previously, the pressure is condensed at the element level, and consequently the variation $Dp\cdot\Delta\mathbf{H}$ can be computed using \eqref{eqn:3.34}. Substituting \eqref{eqn:3.36} into \eqref{eqn:3.35} and using \eqref{eqn:3.9} leads to the set of linearized equations \eqref{eqn:3.10}, with \eqref{eqn:3.10}$_{2}$ replaced by the following equation:
\begin{equation}\label{eqn:3.37}
\begin{aligned}
&\int_{\mathcal{B}}\delta\mathbf{h}^{s}:\left[\mathbb{C}^{d}+pJ\left(\mathbf{I}\otimes\mathbf{I}-2\mathcal{I}\right)\right]:\Delta\mathbf{h}^{s}dV+
\int_{\mathcal{B}}\left(\Delta\mathbf{h}\,\bm{\widehat{\tau}}\right):\delta\mathbf{h}\,dV+
\int_{\mathcal{B}}\left(Dp\cdot\Delta\mathbf{H}\right)J\left(\mathrm{tr}\:\delta\mathbf{h}^{s}\right)\,dV \\
& \quad -\int_{\mathcal{B}}\Delta\bm{\tau}:\delta\mathbf{h}\,dV=
 -\left(\int_{\mathcal{B}}\bm{\widehat{\tau}}:\delta\mathbf{h}^{s}\,dV-
\int_{\mathcal{B}}\bm{\tau}:\delta\mathbf{h}\,dV\right).
\end{aligned}
\end{equation}

Following \citet{SIM91}, we discretize $p$ and $\delta p$ using the relations
\begin{equation}\label{eqn:3.38}
p=\bm{\Gamma}^{\mathsf{T}}_{p}\mathbf{p}\,,\qquad
\delta p=\bm{\Gamma}^{\mathsf{T}}_{p}\delta\mathbf{p}\,,
\end{equation}
where the vectors $\mathbf{p}$ and $\delta\mathbf{p}$ list the degrees of freedom corresponding to $p$ and $\delta p$, respectively and $\bm{\Gamma}_{p}$ is the vector of shape functions that interpolate the pressure and dilation. Similarly, we consider the following forms for discretizing $\theta$ and $\delta\theta$:
\begin{equation}\label{eqn:3.39}
\theta=\bm{\Gamma}^{\mathsf{T}}_{p}\mathbf{q}\,,\qquad
\delta\theta=\bm{\Gamma}^{\mathsf{T}}_{p}\delta\mathbf{q}\,.
\end{equation}
In these equations, the vectors $\mathbf{q}$ and $\delta\mathbf{q}$ contain the degrees of freedom corresponding to $\theta$ and $\delta\theta$. Substituting \eqref{eqn:3.38} and \eqref{eqn:3.39} into \eqref{eqn:3.34} and noting that the vectors $\delta\mathbf{p}$ and $\delta\mathbf{q}$ are arbitrary lead to the following equations for computing the pressure $p$ and dilation $\theta$ at a given point $P\left(\zeta,\eta\right)\in\widehat{\mathscr{Q}}$:
\begin{equation}\label{eqn:3.40}
p=\bm{\Gamma}^{\mathsf{T}}_{p}\mathbf{H}^{-1}_{p}\left[\int_{\mathcal{B}}\bm{\Gamma}U'\left(\theta\right)dV\right]\,,\qquad
\theta=\bm{\Gamma}^{\mathsf{T}}_{p}\mathbf{H}^{-1}_{p}\left(\int_{\mathcal{B}}\bm{\Gamma}JdV\right)\,,
\end{equation}
where
\begin{equation}\label{eqn:3.41}
\mathbf{H}_{p}=\int_{\mathcal{B}}\bm{\Gamma}_{p}\bm{\Gamma}^{\mathsf{T}}_{p}dV.
\end{equation}
From \eqref{eqn:3.40}$_{2}$, it is clear that $\theta$ has a variation with respect to $\Delta\mathbf{H}$. Observing that $DJ\cdot\Delta\mathbf{H}=J\left(\mathrm{tr}\:\Delta\mathbf{h}\right)$ we can express $D\theta\cdot\Delta\mathbf{H}$ in the following form:
\begin{equation}\label{eqn:3.42}
D\theta\cdot\Delta\mathbf{H}=\bm{\Gamma}^{\mathsf{T}}_{p}\mathbf{H}^{-1}_{p}\left[\int_{\mathcal{B}}\bm{\Gamma}J\left(\mathrm{tr}\:\Delta\mathbf{h}\right)dV\right]=\overline{\mathrm{tr}}\:\Delta\mathbf{h}\,,
\end{equation}
where $\overline{\mathrm{tr}}$ is the discrete trace operator. We can now compute $Dp\cdot\Delta\mathbf{H}$ with the help of \eqref{eqn:3.42} as
\begin{equation}\label{eqn:3.43}
Dp\cdot\Delta\mathbf{H}=\bm{\Gamma}^{\mathsf{T}}_{p}\mathbf{H}^{-1}_{p}\left[\int_{\mathcal{B}}\bm{\Gamma}U''\left(\theta\right)\left(\overline{\mathrm{tr}}\:\Delta\mathbf{h}\right)dV\right]\,.
\end{equation}
Substituting \eqref{eqn:3.43} into the third integral on the left-hand side of \eqref{eqn:3.37} gives us
\begin{equation}\label{eqn:3.44}
\int_{\mathcal{B}}\left(Dp\cdot\Delta\mathbf{H}\right)J\left(\mathrm{tr}\:\delta\mathbf{h}^{s}\right)\,dV=
\int_{\mathcal{B}}U''\left(\theta\right)\left(\overline{\mathrm{tr}}\:\delta\mathbf{h}\right)\left(\overline{\mathrm{tr}}\:\Delta\mathbf{h}\right)dV\,.
\end{equation}

For a quadrilateral element $\mathscr{Q}$ with a bilinear displacement field, we assume the pressure $p$ and dilation $\theta$ to be constant throughout the element. Consequently, the shape-function vector used to interpolate the pressure and dilation reduces to $\bm{\Gamma}=\bm 1$. From \eqref{eqn:3.41}, we also observe that the matrix $\mathbf{H}_{p}$ reduces to the material volume $V^{\mathscr{Q}}$ of the element. With these simplifications, \eqref{eqn:3.40} reduces to
\begin{equation}\label{eqn:3.45}
p=U'\left(\theta\right),\qquad\theta=\frac{v^{\mathscr{Q}}}{V^{\mathscr{Q}}}\,,
\end{equation}
where $v^{\mathscr{Q}}$ is the volume of the element in the current configuration. The discretized Eqs.~\eqref{eqn:3.22} and \eqref{eqn:3.26}--\eqref{eqn:3.31} remain valid for incompressible solids. However, Eqs.~\eqref{eqn:3.37} and \eqref{eqn:3.44} indicate that the contribution of pressure must be incorporated into the stiffness matrix $\mathbf{K}^{\alpha\alpha}_{m}$ in \eqref{eqn:3.23}, as demonstrated by the following equation:
\begin{equation}\label{eqn:3.46}
\mathbf{K}^{\alpha\alpha}_{m}=
\int_{\mathcal{B}}\mathbf{\bar{g}}^{c\,\mathsf{T}}\mathbf{D}\,\mathbf{\bar{g}}^{c}dV+
\int_{\mathcal{B}}\mathbf{\bar{g}}^{c\,\mathsf{T}}\mathbf{\widehat{T}}\,\mathbf{\bar{g}}^{c}dV+
\int_{\mathcal{B}}U''\left(\theta\right)\mathbb{\bar{V}}\mathbb{\bar{V}}^{\mathsf{T}}dV\,,
\end{equation}
where from \eqref{eqn:3.42} we have
\begin{equation}\label{eqn:3.47}
\mathbb{\bar{V}}=\frac{1}{V^{\mathscr{Q}}}\int_{\mathcal{B}}J\mathbb{V}dV\,,
\end{equation}
and the vector $\mathbb{V}$ is defined as follows:
\begin{equation}\label{eqn:3.48}
\mathbb{V}^{\mathsf{T}}=
\begin{bmatrix}
  \mathbf{V}^{e_{1}\mathsf{T}}_{1} & \mathbf{V}^{e_{1}\mathsf{T}}_{1} & \mathbf{V}^{e_{1}\mathsf{T}}_{2} & \mathbf{V}^{e_{1}\mathsf{T}}_{2} & \dots & \mathbf{V}^{\mathscr{Q}\mathsf{T}}_{4} & \mathbf{V}^{\mathscr{Q}\mathsf{T}}_{4}
\end{bmatrix}_{1\times24}\,.
\end{equation}
In \eqref{eqn:3.48}, $\mathbf{V}^{e_{i}}_{j},~i=1,\dots,4,~j=1,2$ and $\mathbf{V}^{\mathscr{Q}}_{j},~j=1,\dots,4$ are the global shape functions defined in \eqref{eqn:2.11} and \eqref{eqn:2.13}. We also note that the matrix $\mathbf{D}$ in the first integral of \eqref{eqn:3.46} is the matrix form of the elasticity tensor $\mathbb{C}^{d}+pJ\left(\mathbf{I}\otimes\mathbf{I}-2\mathcal{I}\right)$ (see \eqref{eqn:3.37}).

\subsection{Solvability of the finite element method}
\label{sec:3.4}

Mixed finite element formulations, such as the two formulations discussed in the previous sections, often lead to a system of equations of the form \eqref{eqn:3.28}, where the coefficient matrix has the structure given in \eqref{eqn:3.29}. Problems with this structure belong to the class of saddle-point problems; the existence and uniqueness of solutions for this class of problems are discussed in detail in \citep{BOF13}. Rearranging \eqref{eqn:3.22}, one observes that the tangent stiffness matrix can be written in the following form:
\begin{equation}\label{eqn:3.49}
\mathbf{K}^{\mathscr{Q}}_{t}=
\begin{bmatrix}
  \mathbf{A} & \mathbf{B}^{\mathsf{T}} \\
  \mathbf{B} & \mathbf{0}_{16\times16}
\end{bmatrix}\,,
\end{equation}
where the matrices $\mathbf{A}$, $\mathbf{B}$, and $\mathbf{B}^{\mathsf{T}}$ are given by
\begin{equation}\label{eqn:3.50}
\mathbf{A}=
\begin{bmatrix}
  \mathbf{K}^{\alpha\alpha}_{m} & \mathbf{0}_{24\times8} \\
  \mathbf{0}_{8\times24}        & \mathbf{0}_{8 \times8}
\end{bmatrix}\,,\qquad
\mathbf{B}=
\begin{bmatrix}
  -\mathbf{K}^{\gamma\alpha}_{i} & \mathbf{K}^{\gamma u}_{i}
\end{bmatrix}\,,\qquad
\mathbf{B}^{\mathsf{T}}=
\begin{bmatrix}
  -\mathbf{K}^{\alpha\gamma}_{i} \\
  \mathbf{K}^{u\gamma}_{i}
\end{bmatrix}\,.
\end{equation}
For a given quadrilateral element, the numbers of degrees of freedom associated with the displacement gradient, stress tensor, and displacement field are denoted by $r=24$, $s=16$, and $k=8$, respectively (see Table~\ref{tab:3.1}). Consequently, the matrices $\mathbf{K}^{\alpha\alpha}_{m}$, $\mathbf{K}^{\gamma\alpha}_{i}$, and $\mathbf{K}^{\gamma u}_{i}$ have dimensions $r\times r$, $s\times r$, and $s\times k$, respectively. It is straightforward to verify that the condition $r+k>s>k$, which is necessary for the invertibility of the tangent stiffness matrix $\mathbf{K}^{\mathscr{Q}}_{t}$, is satisfied \cite[Remark 3.2.1]{BOF13}. However, additional conditions must be imposed on the matrices $\mathbf{K}^{\alpha\alpha}_{m}$, $\mathbf{K}^{\gamma\alpha}_{i}$, and $\mathbf{K}^{\gamma u}_{i}$ in order to ensure that the tangent stiffness matrix is non-singular. These conditions may be summarized as follows \cite[Section 3.2.5]{BOF13}:
\begin{align}\label{eqn:3.51}
\ker\left(\mathbf{K}^{\gamma u}_{i}\right)=\{\mathbf{0}_{k}\}\,,\qquad
\ker\left(\mathbf{K}^{\alpha\gamma}_{i}\right)\cap\ker\left(\mathbf{K}^{u\gamma}_{i}\right)=\{\mathbf{0}_{s}\}\,,\qquad
\ker\left(\mathbf{K}^{\alpha\alpha}_{m}\right)\cap K=\{\mathbf{0}_{r}\}\,,
\end{align}
where the set $K$ is defined by
\begin{equation}\label{eqn:3.52}
K=\Big\{\mathbf{x}\in\mathbb{R}^{r}\;\text{such that}\;\mathbf{y}^{\mathsf{T}}\mathbf{K}^{\gamma\alpha}_{i}\mathbf{x}=0\,,\quad\forall~\mathbf{y}\in\ker\left(\mathbf{K}^{u\gamma}_{i}\right)\Big\}\,.
\end{equation}
The conditions in \eqref{eqn:3.51} are necessary and sufficient for the tangent stiffness matrix $\mathbf{K}^{\mathscr{Q}}_{t}$ to be non-singular \citep{BOF13}.

Considering a single square element in the physical domain, one observes that $\dim\!\left(\operatorname{im}\!\left(\mathbf{K}^{\gamma u}_{i}\right)\right)=6<k=8$. However, after imposing the boundary conditions, one obtains $\dim\!\left(\operatorname{im}\!\left(\mathbf{K}^{\gamma u}_{i}\right)\right)=k=5$, and hence condition \eqref{eqn:3.51}$_1$ is satisfied.

\begin{remark}\label{rem:3.1}
We note that the satisfaction of condition \eqref{eqn:3.51}$_{1}$ is due to the fact that all $16$ stress degrees of freedom are defined on the edges of the quadrilateral. Alternatively, one could use the $\mathcal{ABF}_{0}$ shape functions for the stress tensor, with $8$ degrees of freedom on the edges and $4$ degrees of freedom on the quadrilateral itself \citep{ARN05}. However, in this case, $\dim\left(\operatorname{im}\left(\mathbf{K}^{\gamma u}_{i}\right)\right)=4$, and condition \eqref{eqn:3.51}$_{1}$ cannot be satisfied even after imposing the boundary conditions. In other words, the form of the shape functions used to interpolate the stress tensor and the number of stress degrees of freedom assigned to the edges are chosen so that condition \eqref{eqn:3.51}$_{1}$ is satisfied.
\end{remark}

It is also observed that $\dim\left(\operatorname{im}\left(\mathbf{K}^{\alpha\gamma}_{i}\right)\right)=s=16$ (equivalently, $\ker\left(\mathbf{K}^{\alpha\gamma}_{i}\right)=\{\mathbf{0}_{s}\}$), and therefore condition \eqref{eqn:3.51}$_{2}$ is satisfied. Regarding the set $K$ in \eqref{eqn:3.52}, we note that the space $\mathbb{R}^{s}$ admits the following direct-sum decomposition:
\begin{equation}\label{eqn:3.53}
\mathbb{R}^{s}=\operatorname{im}\left(\mathbf{K}^{\gamma u}_{i}\right)\oplus
\left(\operatorname{im}\left(\mathbf{K}^{\gamma u}_{i}\right)\right)^{\bot}=
\operatorname{im}\left(\mathbf{K}^{\gamma u}_{i}\right)\oplus
\ker\left(\mathbf{K}^{u\gamma}_{i}\right).
\end{equation}
For every vector $\mathbf{y}\in\mathbb{R}^{s}$, we can define its orthogonal projection onto $\operatorname{im}\left(\mathbf{K}^{\gamma u}_{i}\right)$ by $\mathbf{P}_{\gamma u}\mathbf{y}$, where the projection operator $\mathbf{P}_{\gamma u}$ is defined as
\begin{equation}\label{eqn:3.54}
\mathbf{P}_{\gamma u}=\mathbf{K}^{\gamma u}_{i}\mathbf{K}^{\gamma u\,\dag}_{i},
\end{equation}
where $\mathbf{K}^{\gamma u\,\dag}_{i}=\left(\mathbf{K}^{\gamma u\,\mathsf{T}}_{i}\mathbf{K}^{\gamma u}_{i}\right)^{-1}\mathbf{K}^{\gamma u\,\mathsf{T}}_{i}$ is the Moore-Penrose inverse of $\mathbf{K}^{\gamma u}_{i}$ \citep{ODE10}. On the other hand, $\left(\mathbf{I}-\mathbf{P}_{\gamma u}\right)\mathbf{y}$ projects the vector $\mathbf{y}$ onto $\ker\left(\mathbf{K}^{u\gamma}_{i}\right)$. Based on these definitions, we can define the set $K$ in the alternative form as:
\begin{equation}\label{eqn:3.55}
K=\Big\{\mathbf{x}\in\mathbb{R}^{r}\;\text{such that}\;\left(\mathbf{I}-\mathbf{P}_{\gamma u}\right)\mathbf{K}^{\gamma\alpha}_{i}\mathbf{x}=0\Big\}\,.
\end{equation}
Both definitions of the set $K$ in \eqref{eqn:3.52} and \eqref{eqn:3.55} imply that the orthogonal projection of $\mathbf{K}^{\gamma\alpha}_{i}\mathbf{x}$ onto $\ker\left(\mathbf{K}^{u\gamma}_{i}\right)$ must vanish. Observing \eqref{eqn:3.53}, this in turn implies that $\mathbf{K}^{\gamma\alpha}_{i}\mathbf{x}\in\operatorname{im}\left(\mathbf{K}^{\gamma u}_{i}\right)$. To obtain a basis for $K$, one can compute the null space of the matrix $\left(\mathbf{I}-\mathbf{P}_{\gamma u}\right)\mathbf{K}^{\gamma\alpha}_{i}$. Condition \eqref{eqn:3.51}$_{3}$ is equivalent to $\mathbf{x}^{\mathsf{T}}_{j}\mathbf{K}^{\alpha\alpha}_{m}\mathbf{x}_{j}>0,~j=1,\dots,s-k$ for every vector $\mathbf{x}_{j}$ belonging to this basis. Numerical experiments using the matrix $\mathbf{K}^{\alpha\alpha}_{m}$ defined in \eqref{eqn:3.23}$_{2}$ indicate that this condition, and consequently condition \eqref{eqn:3.51}$_{3}$, is satisfied. The satisfaction of the necessary and sufficient conditions in \eqref{eqn:3.51} ensures the invertibility of the tangent stiffness matrix $\mathbf{K}^\mathscr{Q}_{t}$ in \eqref{eqn:3.29}. In fact, a singular value decomposition of $\mathbf{K}^\mathscr{Q}_{t}$ shows that this matrix is positive definite after the imposition of the boundary conditions. Finally, it is important to note that the conditions in \eqref{eqn:3.51} can be expressed in the form of inf--sup, or Ladyzhenskaya--Babu\v{s}ka--Brezzi, conditions \citep{SHO18}. Therefore, the element developed in this work satisfies these conditions as well.

\section{Numerical Examples}
\label{sec:4}

In this section, we present several numerical examples to study the stability of the compressible and incompressible formulations of the first-order CSMFE quadrilateral element. The performance of the element, including the deformation and load-displacement curves it predicts, is compared with the CMSFE formulations developed for simplicial elements \citep{JAH22,JAH25} and with the classical displacement-based and U/P mixed quadrilateral elements \citep{BAT96,ZIE05,HUG00,BEL05,WRI08}. The problems presented in this section have been widely used as benchmark problems for assessing the numerical stability of finite elements. The solution of some of these problems may lead to the emergence of spurious energy modes in first-order elements. In other problems, shear locking and numerical instabilities may arise in the near-incompressible regime. The numerical results obtained using the first-order quadrilateral element developed in this work demonstrate its excellent performance and the absence of such numerical instabilities. All numerical examples are solved using the full Newton-Raphson method together with the tangent stiffness matrix derived in \eqref{eqn:3.29}. A convergence tolerance of $10^{-9}$ is used in all numerical examples. Different load increments are employed in the simulations, and the numerical experiments show that, on average, $3$ equilibrium iterations are required for convergence at each load step in the compressible formulation. However, this number increases to $6$ iterations for the incompressible formulation.

\begin{remark}\label{rem:4.1}
It is important to note that the first-order $U/P$ mixed formulation does not produce correct results for certain numerical examples presented below. It was therefore necessary to compare the results of the incompressible formulation with those obtained using the second-order $U/P$ mixed formulation. It is observed that the results obtained with the first-order quadrilateral element developed in this paper agree well with those of the second-order $U/P$ mixed formulation.
\end{remark}

Two material models are used in the numerical examples. The first model is employed in the compressible formulation, whereas the second model is used in the incompressible formulation. These materials are characterized by the following strain energy functions:
\begin{equation}\label{eqn:4.1}
\begin{aligned}
\widehat{W}_{1}\left(\mathbf{C}\right)&=\frac{\mu}{2}\Bigl[\left(I_{1}-3\right)-\ln J^{2}\Bigr]+
\frac{\kappa}{2}\left(J-1\right)^{2}\,,\\
\widehat{W}_{2}\left(\mathbf{\bar{C}}\right)&=\frac{\mu}{2}\left(\bar{I}_{1}-3\right)\,,\qquad
U\left(\theta\right)=\frac{\kappa}{2}\left(\theta-1\right)^{2}\,,
\end{aligned}
\end{equation}
where $\mu$ is the shear modulus, $\kappa$ is the bulk modulus, and $I_{1}=\mathrm{tr}\:\mathbf{C}$ and $\bar{I}_{1}=\mathrm{tr}\:\mathbf{\bar{C}}$ are the first invariants of $\mathbf{C}$ and $\mathbf{\bar{C}}$, respectively. The Kirchhoff stress tensors corresponding to the two material models have the following representations:
\begin{equation}\label{eqn:4.2}
\begin{aligned}
\widehat{\bm{\tau}}_{1}&=\mu\mathbf{B}+
\Bigl[\kappa J\left(J-1\right)-\mu\Bigr]\mathbf{I}\,,\\
\widehat{\bm{\tau}}_{2}=J^{-\frac{2}{3}}\left(\operatorname{dev}\:\bm{\bar{\tau}}\right)&+Jp\mathbf{I}\,,\qquad
\bm{\bar{\tau}}=\mu\mathbf{\bar{B}}\,,\qquad p=\kappa\left(\theta-1\right)\,,
\end{aligned}
\end{equation}
where $\mathbf{B}=\mathbf{F}\mathbf{F}^{\mathsf{T}}$ is the left Cauchy-Green deformation tensor and $\mathbf{\bar{B}}=J^{-\frac{2}{3}}\mathbf{B}$. We note that $\bm{\widehat{\tau}}_{2}$ is the Kirchhoff stress tensor obtained in \eqref{eqn:3.33} from the variation of $\Pi_{2}$ with respect to $\delta\mathbf{H}$ (see the discussion preceding \eqref{eqn:3.33}). It should also be noted that $\bm{\widehat{\tau}}_{1}$ and $\bm{\widehat{\tau}}_{2}$ are the constitutively derived Kirchhoff stress tensors. This is in contrast to the interpolated Kirchhoff stress tensor $\bm{\tau}$ obtained using \eqref{eqn:3.11}$_{2}$. The spatial elasticity tensors corresponding to the Kirchhoff stress tensors in \eqref{eqn:4.2} are
\begin{equation}\label{eqn:4.3}
\mathbb{C}_{1}=\kappa J\left(2J-1\right)\mathbf{I}\otimes\mathbf{I}+
2\Bigl[\mu-\kappa J\left(J-1\right)\Bigr]\boldsymbol{\mathcal{I}}\,,\qquad
\mathbb{C}_{2}=\mathbb{C}^{d}+pJ\left(\mathbf{I}\otimes\mathbf{I}-2\mathcal{I}\right)\,.
\end{equation}
The elasticity tensor $\mathbb{C}_{2}$ in \eqref{eqn:4.3}$_{2}$ is the same tensor that appears in \eqref{eqn:3.36}$_{2}$. The tensor $\mathbb{C}^{d}$ can be derived from \eqref{eqn:4.1}$_{2}$ as follows:
\begin{equation}\label{eqn:4.4}
\mathbb{C}^{d}=-\frac{2}{3}\mu J^{-\frac{2}{3}}\Big[\left(\mathrm{dev}\:\mathbf{B}\right)\otimes\mathbf{I}+\mathbf{I}\otimes\left(\mathrm{dev}\:\mathbf{B}\right)\Big]
-\frac{2}{9}\mu J^{-\frac{2}{3}}\left(\mathrm{tr}\:\mathbf{B}\right)\mathbf{I}\otimes\mathbf{I}
+\frac{2}{3}\mu J^{-\frac{2}{3}}\left(\mathrm{tr}\:\mathbf{B}\right)\mathcal{I}\,.
\end{equation}

\subsection{Cook's membrane problem}
\label{sec:4.1}

\begin{figure}
\centering
\includegraphics[width=0.3\textwidth]{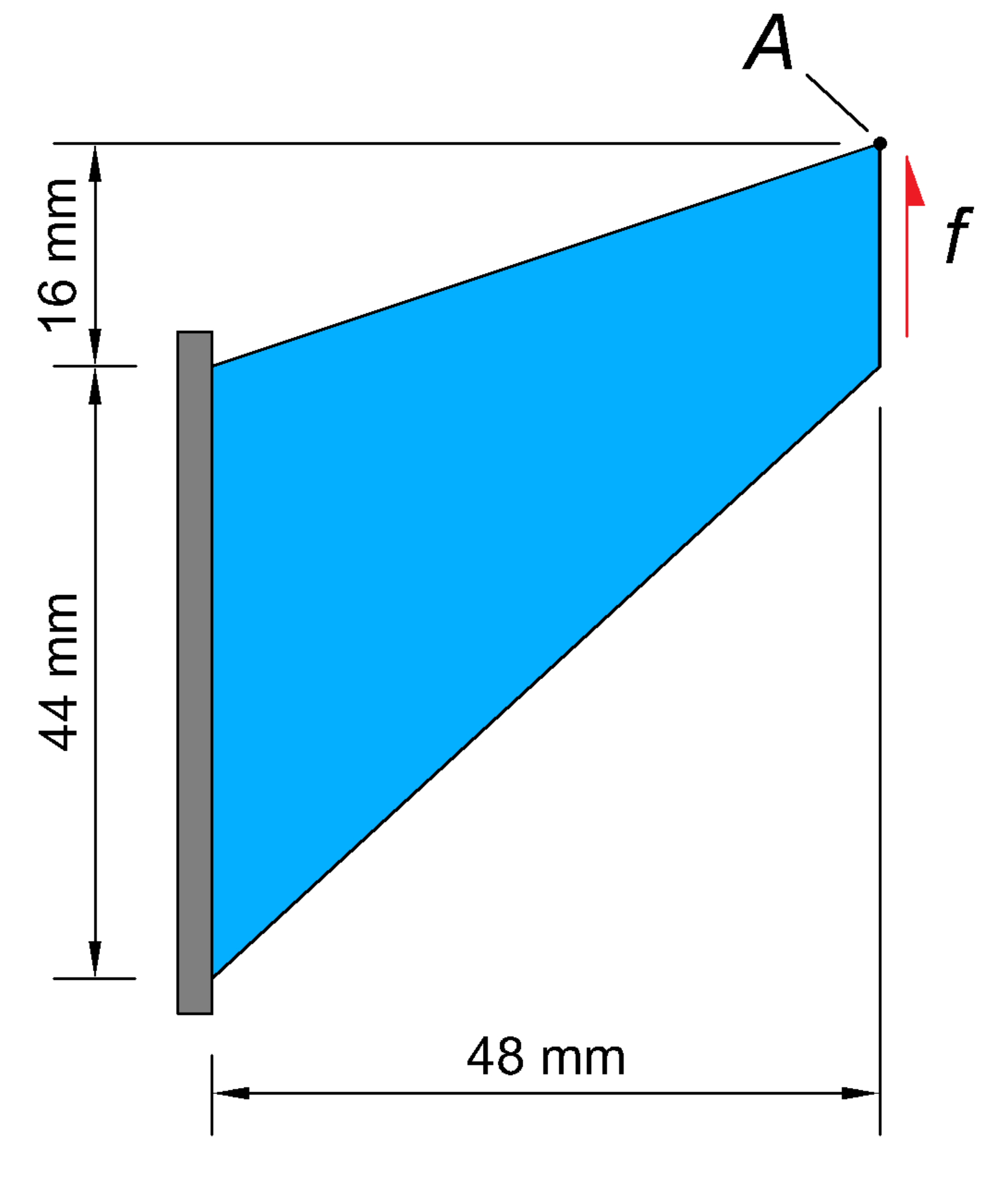}
\vskip 0.1in
\caption{Geometry and boundary conditions for Cook's membrane problem: The left edge is fixed, while the right edge is subject to a uniform shearing load $f$.}
\label{fig:41}
\end{figure}

The classical Cook's membrane problem has been studied extensively in the literature \citep{BRI96,GLA97,REE00,ANG17,DHA22a,JAH22,JAH25}. It is a benchmark problem for assessing the performance of finite elements in bending-dominated plane strain deformations. A tapered panel, fixed at one end, is subjected to a uniform shearing load at the opposite end. The geometry and boundary conditions are shown in Figure~\ref{fig:41}. Numerical experiments indicate that the first CSMFE formulation together with the first material model in \eqref{eqn:4.1} exhibits an overly stiff response and leads to the emergence of a checkerboard pattern in the stress contours. However, the problem is modeled correctly using the second CSMFE formulation and the second material model. The material parameters used in this model are $\mu = 80.194$ MPa and $\kappa = 400,889.8$ MPa. It is known that certain first-order elements, such as Q1 and Q1/P0, as well as some enhanced strain-based elements, such as Q1/ES4, exhibit overly stiff behavior when solving Cook's membrane problem \citep{GLA97,REE00}.

\begin{figure}
\centering
\includegraphics[width=0.65\textwidth]{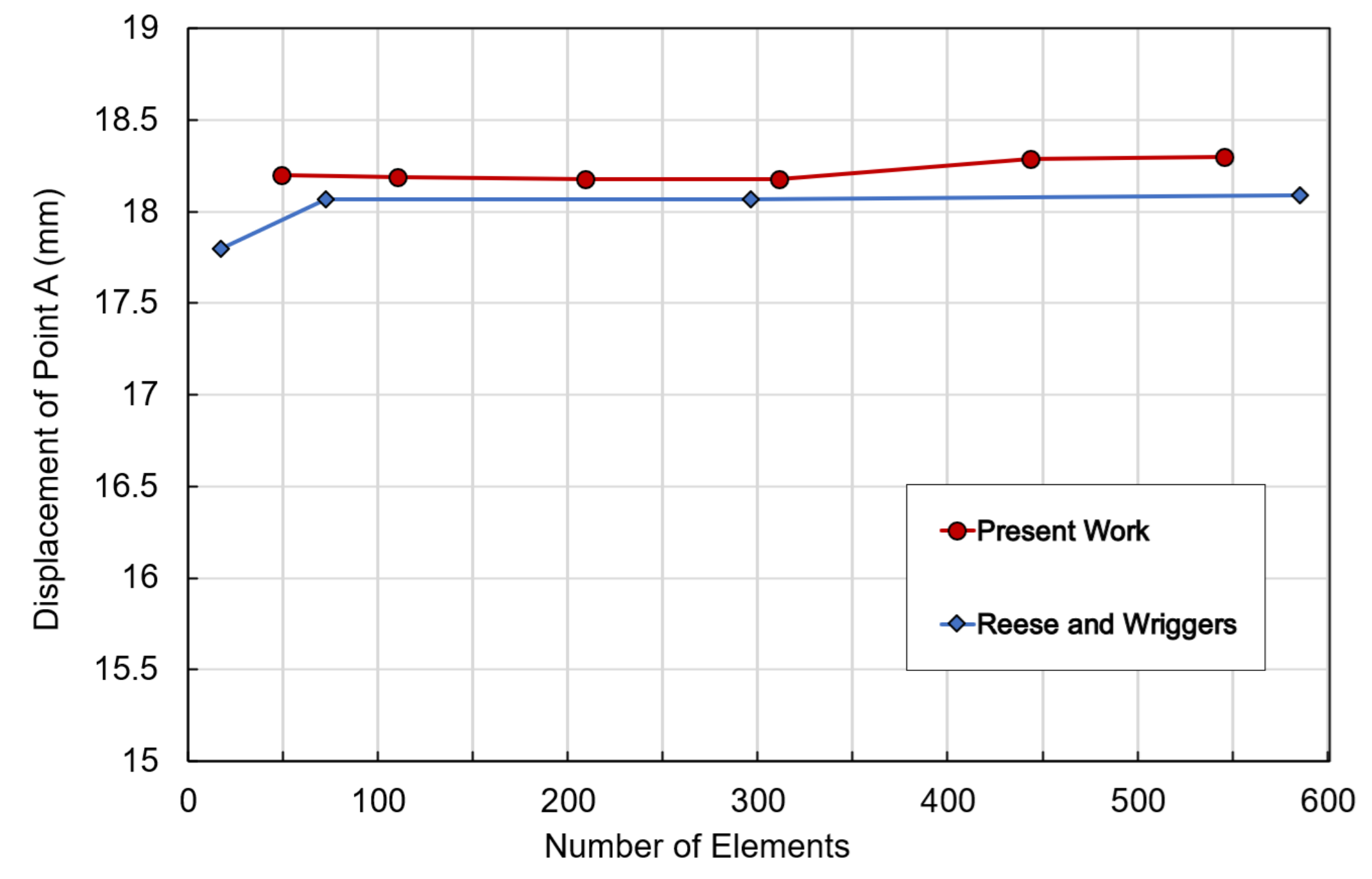}
\vskip 0.1in
\caption{Comparison of the vertical displacement of point $A$ for Cook's membrane problem with \citep{REE00} using different number of elements and $f=24$ MPa.}
\label{fig:42}
\end{figure}

\begin{figure}
\centering
\includegraphics[width=0.8\textwidth]{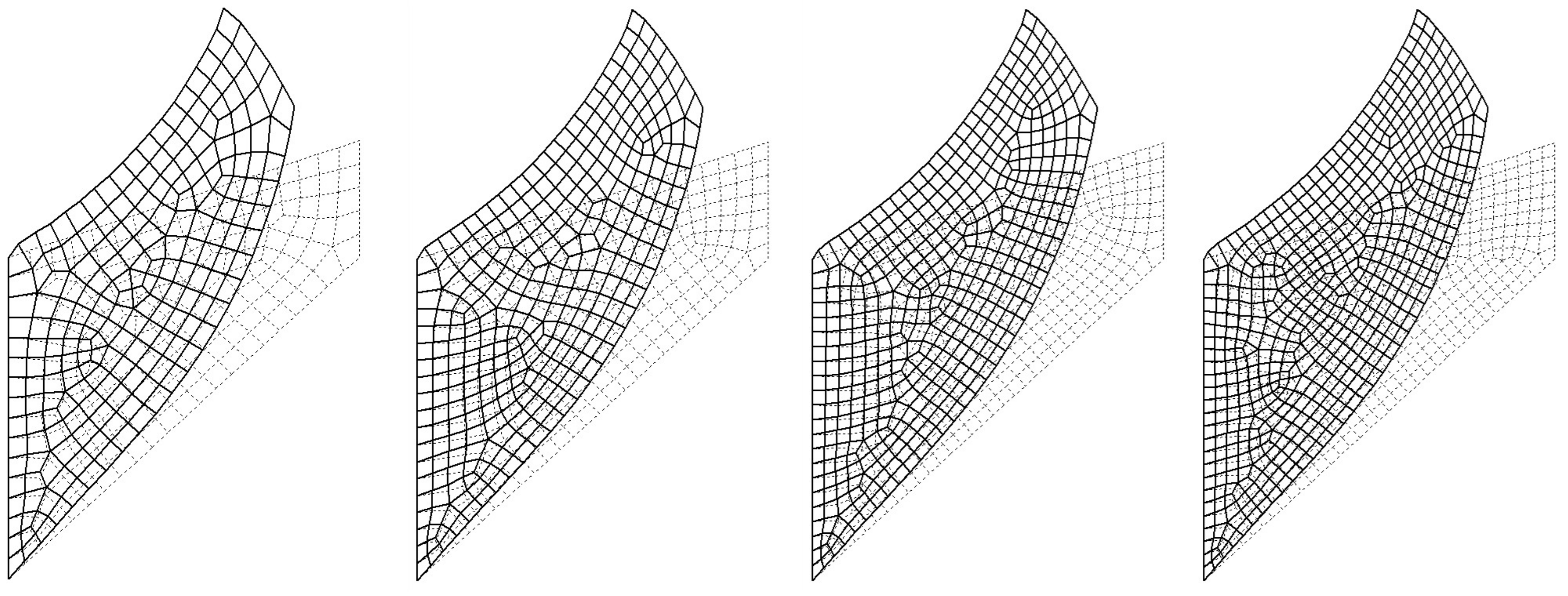}
\vskip 0.1in
\caption{Deformed and undeformed configurations of Cook's membrane problem for meshes $3$--$6$ of the quadrilateral element developed in this work. The shearing load is $f=24$ MPa for all meshes.}
\label{fig:43}
\end{figure}

\begin{figure}
\centering
\includegraphics[width=0.5\textwidth]{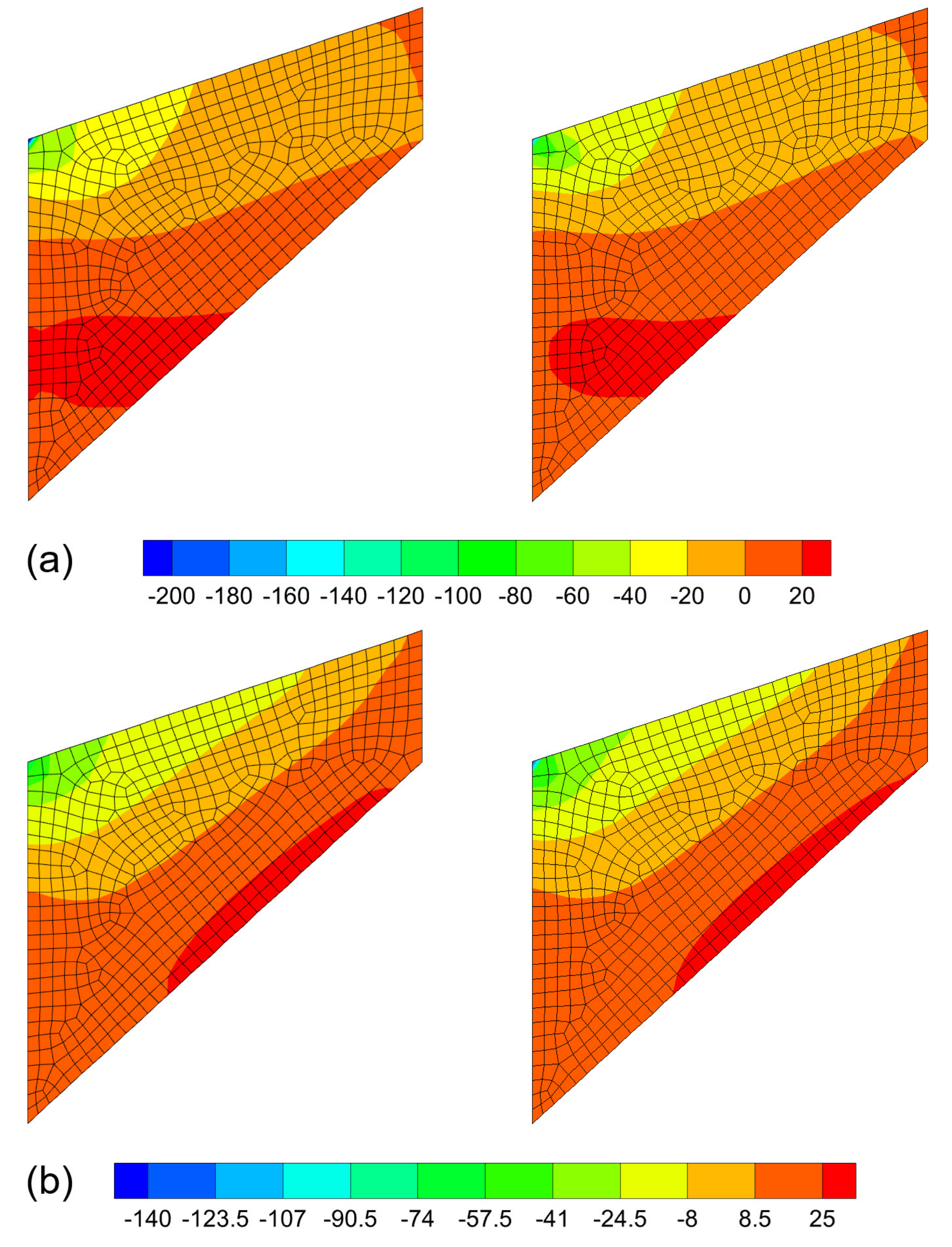}
\vskip 0.1in
\caption{Cook's membrane problem: (a) Contours of the Kirchhoff stress component $\tau_{11}$~(MPa) using the element developed in this work and the second-order mixed U/P finite element. (b) Contours of the pressure $p$~(MPa) using the element developed in this work and the second-order mixed U/P finite element.}
\label{fig:44}
\end{figure}

Several meshes are considered to solve this problem. The number of elements and the corresponding number of degrees of freedom for these meshes are listed in Table~\ref{tab:4.1}. A uniform shearing load of $f=24$ MPa is applied at the free end in $1000$ load steps. The resulting vertical displacement of point $A$ computed using the meshes in Table~\ref{tab:4.1} is compared in Figure~\ref{fig:42} with the results obtained using the Q1SP elements of \citet{REE00}. The corresponding values of the displacement of point $A$ are also reported in Table~\ref{tab:4.1} for the meshes used in the present simulations. The undeformed and deformed configurations for the meshes in Table~\ref{tab:4.1} are shown in Figure~\ref{fig:43}. The contours of the stress component $\tau_{11}$ and the pressure $p$ are shown in Figure~\ref{fig:44} for both the element developed in this work and the second-order mixed U/P finite element. An excellent agreement is observed between the results obtained with the two elements.

\begin{table}
\centering
\caption{The number of elements, the number of degrees of freedom, and the displacement of point $A$ for the meshes of quadrilateral elements used to model Cook's membrane problem. The value of displacement at point $A$ corresponds to a uniform shearing load of $f=24$ MPa.}
\label{tab:4.1}
\renewcommand{\arraystretch}{1.5}
\renewcommand{\tabcolsep}{0.2cm}
\begin{tabular}{c c c c}
\hline
Mesh & No. Elements & DOFs  & Displ. A (mm) \\
\hline
1    & 49           & 1408  & 18.20         \\
2    & 110          & 3058  & 18.19         \\
3    & 209          & 5702  & 18.18         \\
4    & 311          & 8418  & 18.18         \\
5    & 443          & 11914 & 18.29         \\
6    & 545          & 14610 & 18.30         \\
\hline
\end{tabular}
\end{table}

\subsection{Rubber sealing problem}
\label{sec:4.2}

The rubber sealing problem involves the compression of a rubber seal and was originally studied in \citep{BRI98} in the context of a posteriori error estimation for finite element computations in finite elasticity. In that work, the Q2/P1 element together with the Ogden material model was used in the simulations. The applied displacement, however, was limited to $2.0$ mm. The problem is known to be challenging for first-order elements such as Q1 and Q1/P0. Using meshes of Q1/ET4 elements and a compressible neo-Hookean material model, it was possible to simulate the problem up to an applied displacement of $1.8$ mm \citep{JAH22}. First- and second-order simplex CSMFEs have also been applied successfully to model the problem for applied displacements as large as $2.2$ mm \citep{ANG17,JAH22,JAH25}.

\begin{figure}
\centering
\includegraphics[width=0.45\textwidth]{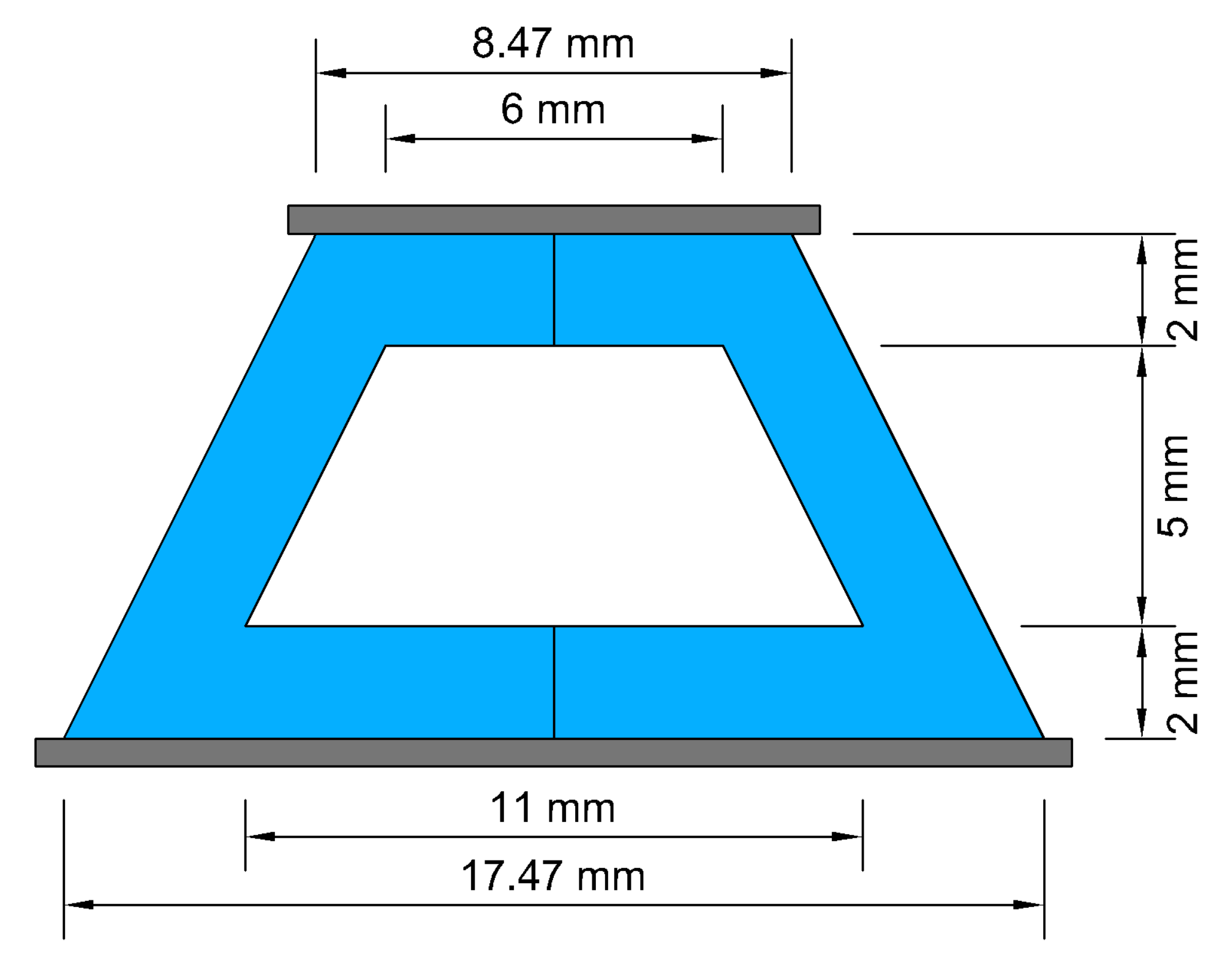}
\vskip 0.1in
\caption{Geometry and boundary conditions for the rubber sealing problem: The bottom edge is fixed, while a specified displacement of $2.2$ mm is imposed on the top edge. Due to the symmetry of the geometry and loading, only the right half of the model is considered.}
\label{fig:45}
\end{figure}

\begin{figure}
\centering
\includegraphics[width=0.5\textwidth]{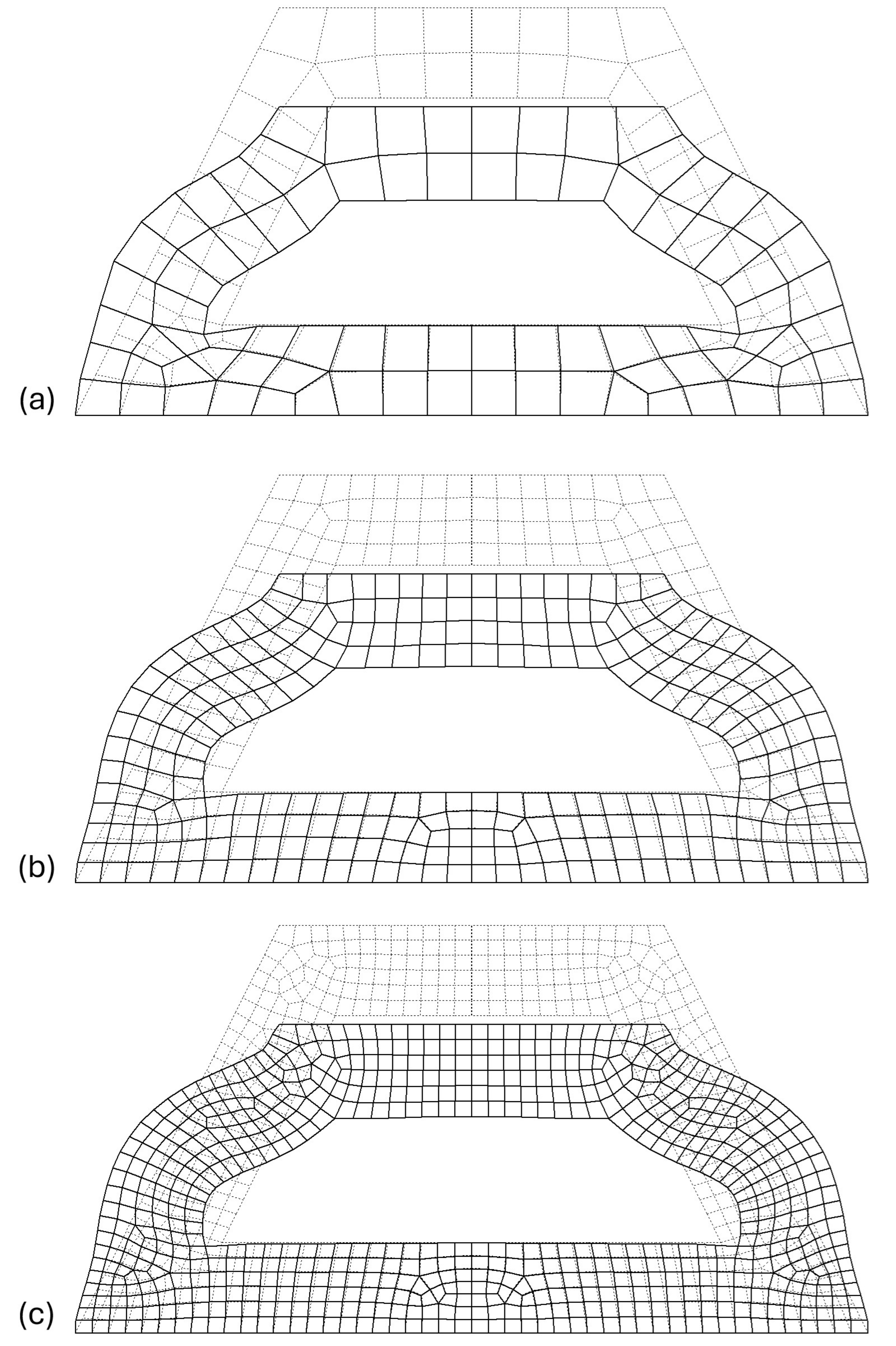}
\vskip 0.1in
\caption{Deformed configurations obtained using the compressible formulation for the meshes used to discretize the rubber sealing problem: (a) $43$, (b) $140$ and (c) $329$ quadrilateral elements.}
\label{fig:46}
\end{figure}

\begin{figure}
\centering
\includegraphics[width=0.70\textwidth]{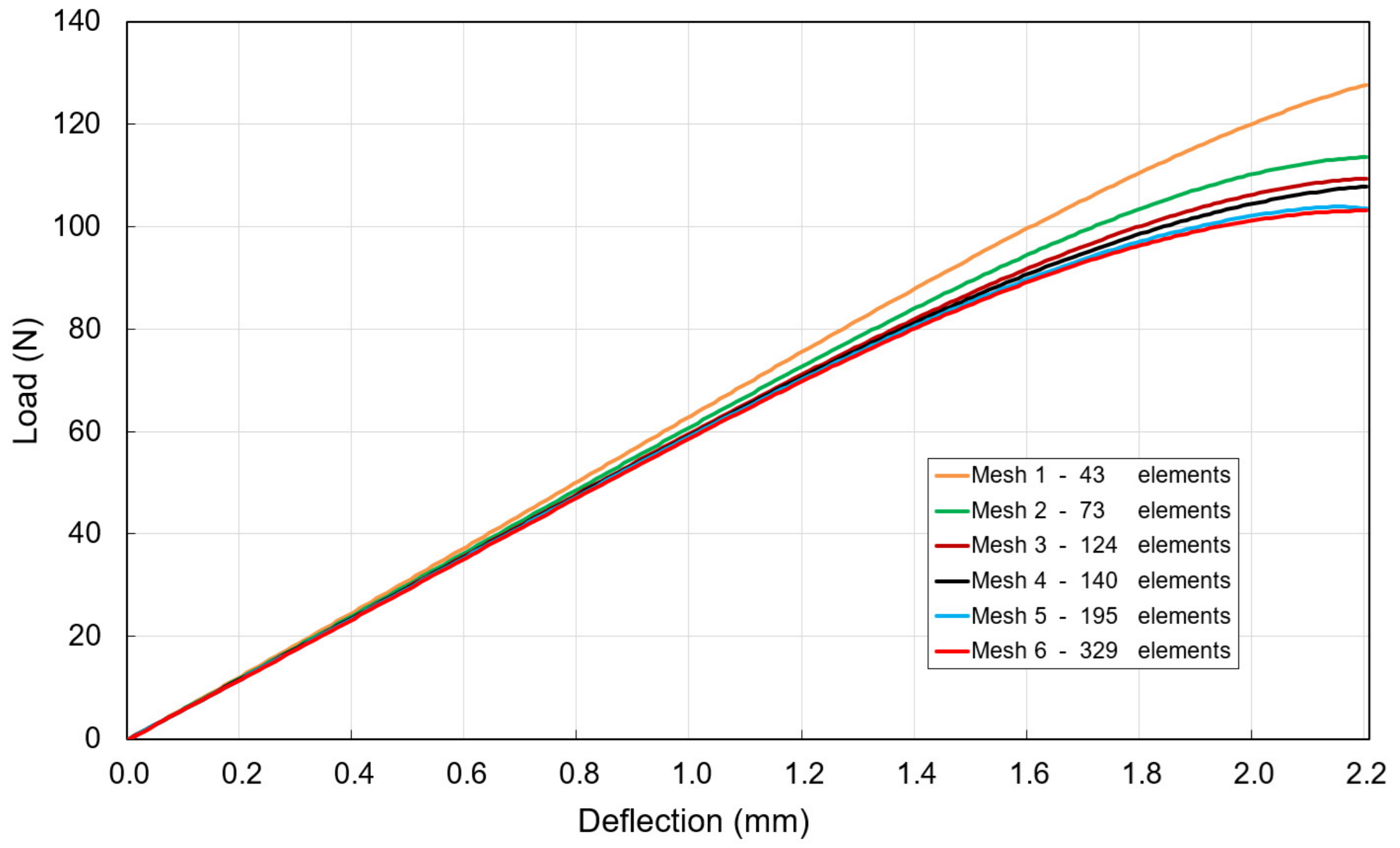}
\vskip 0.1in
\caption{Load-deflection curves for the rubber sealing problem using the compressible formulation, the first material model in \eqref{eqn:4.1}, and different meshes of the quadrilateral element.}
\label{fig:47}
\end{figure}

\begin{figure}
\centering
\includegraphics[width=0.55\textwidth]{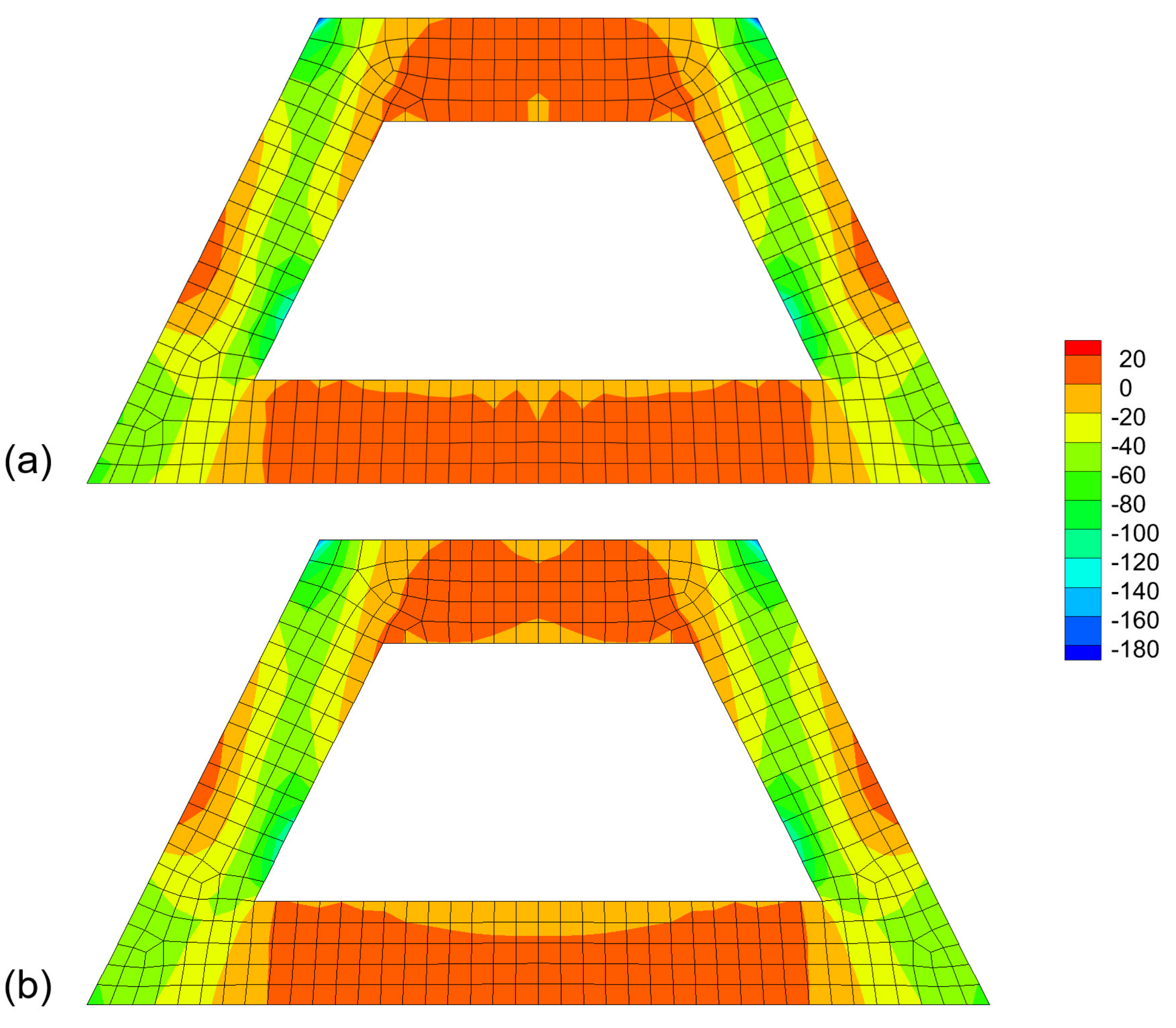}
\caption{Contours of the Kirchhoff stress component $\tau_{22}$ (MPa) for the rubber sealing problem using the incompressible formulation and the second material model in \eqref{eqn:4.1} for a target displacement of $2.0$~mm; (a) $220$ first-order quadrilateral element developed in this work; (b) $220$ second-order mixed U/P finite elements.}
\label{fig:48}
\end{figure}

\begin{figure}
\centering
\includegraphics[width=0.55\textwidth]{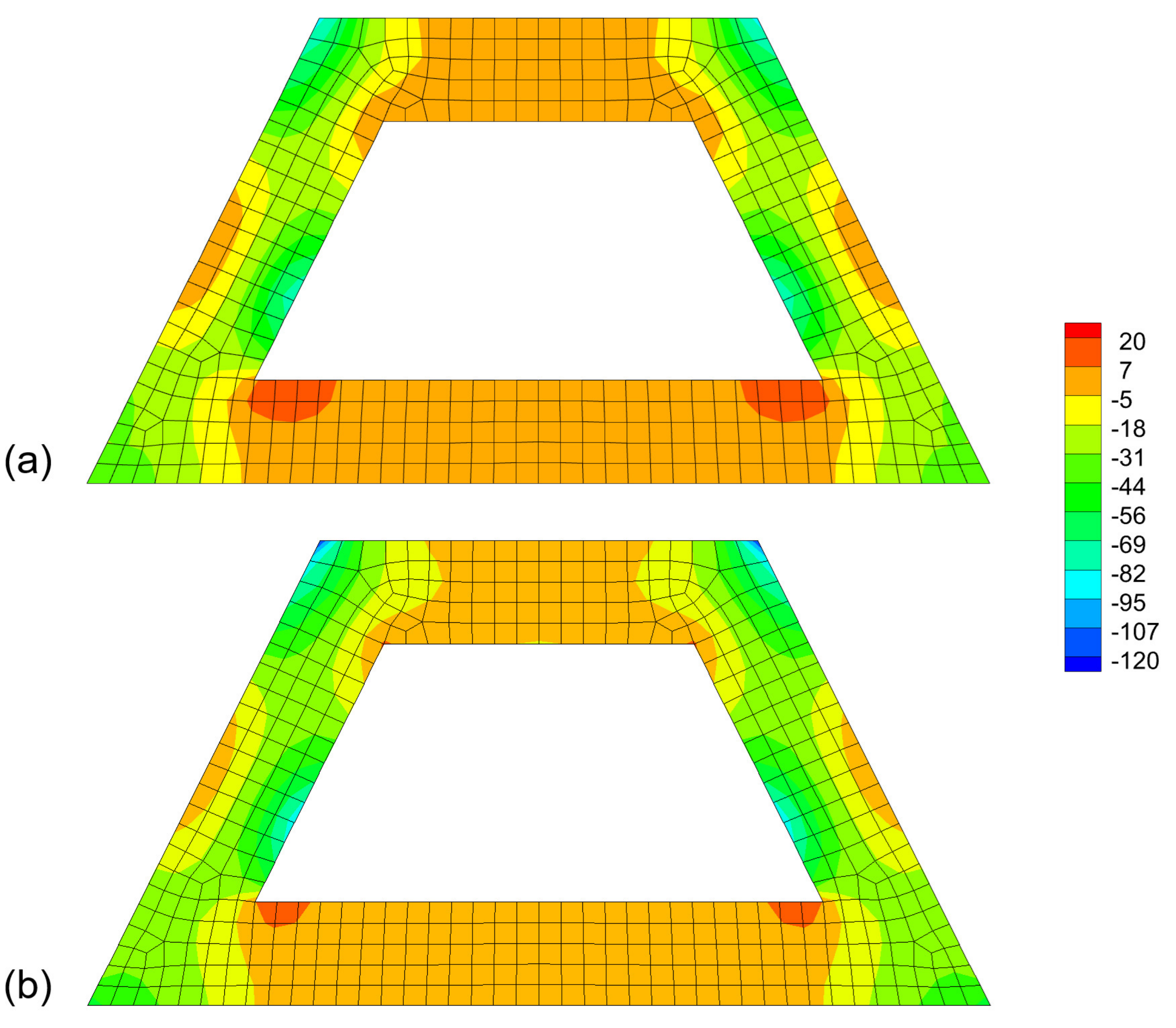}
\caption{Contours of the pressure $p$ (MPa) for the rubber sealing problem using the incompressible formulation and the second material model in \eqref{eqn:4.1} for a target displacement of $2.0$~mm; (a) $220$ first-order quadrilateral element developed in this work; (b) $220$ second-order mixed U/P finite elements.}
\label{fig:49}
\end{figure}

The geometry and boundary conditions of the problem are shown in Figure~\ref{fig:45}. Both the compressible and incompressible formulations, together with the first and second material models in \eqref{eqn:4.1}, are used to solve the problem. The simulations are performed using the material parameters $\mu=82.19$ and $\kappa=1000$ MPa. For the compressible formulation, a downward displacement of $2.2$ mm is prescribed on the top edge, while the bottom edge remains fixed. The prescribed displacement is applied in $220$ load steps. Owing to the symmetry of the geometry and loading, only the right half of the specimen is modeled. Several meshes of the quadrilateral element are considered to study the performance of the proposed element for this problem. The number of elements, the corresponding number of degrees of freedom, and the final value of the load applied to the top edge for these meshes are listed in Table~\ref{tab:4.2}. The undeformed and deformed configurations obtained using the compressible formulation are shown in Figure~\ref{fig:46} for selected meshes. The load-displacement curves corresponding to various meshes and the first material model are plotted in Figure~\ref{fig:47}. The contours of the stress component $\tau_{22}$ and the pressure $p$, obtained using the incompressible formulation and the second material model in \eqref{eqn:4.1}, are shown in Figures~\ref{fig:48} and \ref{fig:49}, respectively, for a target displacement of $2.0$ mm. In these figures, the performance of the first-order quadrilateral element developed in this work is compared with that of the second-order mixed U/P formulation. The figures indicate that there is good agreement between the results obtained using the two elements.

\begin{table}
\centering
\caption{The number of elements, the number of degrees of freedom, and the final value of the load applied to the top edge at the target displacement of $2.2$~mm for the meshes of the quadrilateral elements using the compressible formulation and the first material model.}
\label{tab:4.2}
\renewcommand{\arraystretch}{1.5}
\renewcommand{\tabcolsep}{0.2cm}
\begin{tabular}{c c c c}
\hline
   Mesh &   No. Elements   &   DOFs    &   Load (N)   \\
\hline
   1    &   43             &   1296    &   127.78     \\
   2    &   73             &   2127    &   113.64     \\
   3    &   124            &   3512    &   109.43     \\
   4    &   140            &   3961    &   107.93     \\
   5    &   195            &   5426    &   103.58     \\
   6    &   329            &   9047    &   103.13     \\
\hline
\end{tabular}
\end{table}

\subsection{Tension of a perforated block}
\label{sec:4.4}

This example has been used in various forms to study the behavior of first- and second-order compatible strain simplicial elements under large stretches \citep{SHO18,DHA22a,JAH22,JAH25}. \citet{DHA22a} considered a $2\times2$ mm block subjected to a tensile deformation of $100\%$ in order to investigate the performance of a first-order element. In that work, the material parameters were taken to be $\mu=10$ and $\kappa=1000$ MPa. The edge of the block subjected to the tensile deformation was free to move in the transverse direction. \citet{JAH22} studied a similar problem using a $1\times1$ mm block and the material parameters $\mu=80.192$ and $\kappa=400,933.33$ MPa. As in the work of \citet{DHA22a}, the edge subjected to the tensile deformation was free to move in the transverse direction. The applied deformation, however, was limited to $50\%$ of the height of the block. \citet{JAH25} employed second-order compatible strain simplex elements to model a $1\times1$ mm block subjected to a tensile deformation of $300\%$ applied to its top edge. Unlike the previous studies, the top edge was restrained in the transverse direction. \citet{SHO18} considered the same block geometry and loading conditions as those in \citep{JAH25}, but in the context of incompressible elasticity.

\begin{figure}
\centering
\includegraphics[width=0.3\textwidth]{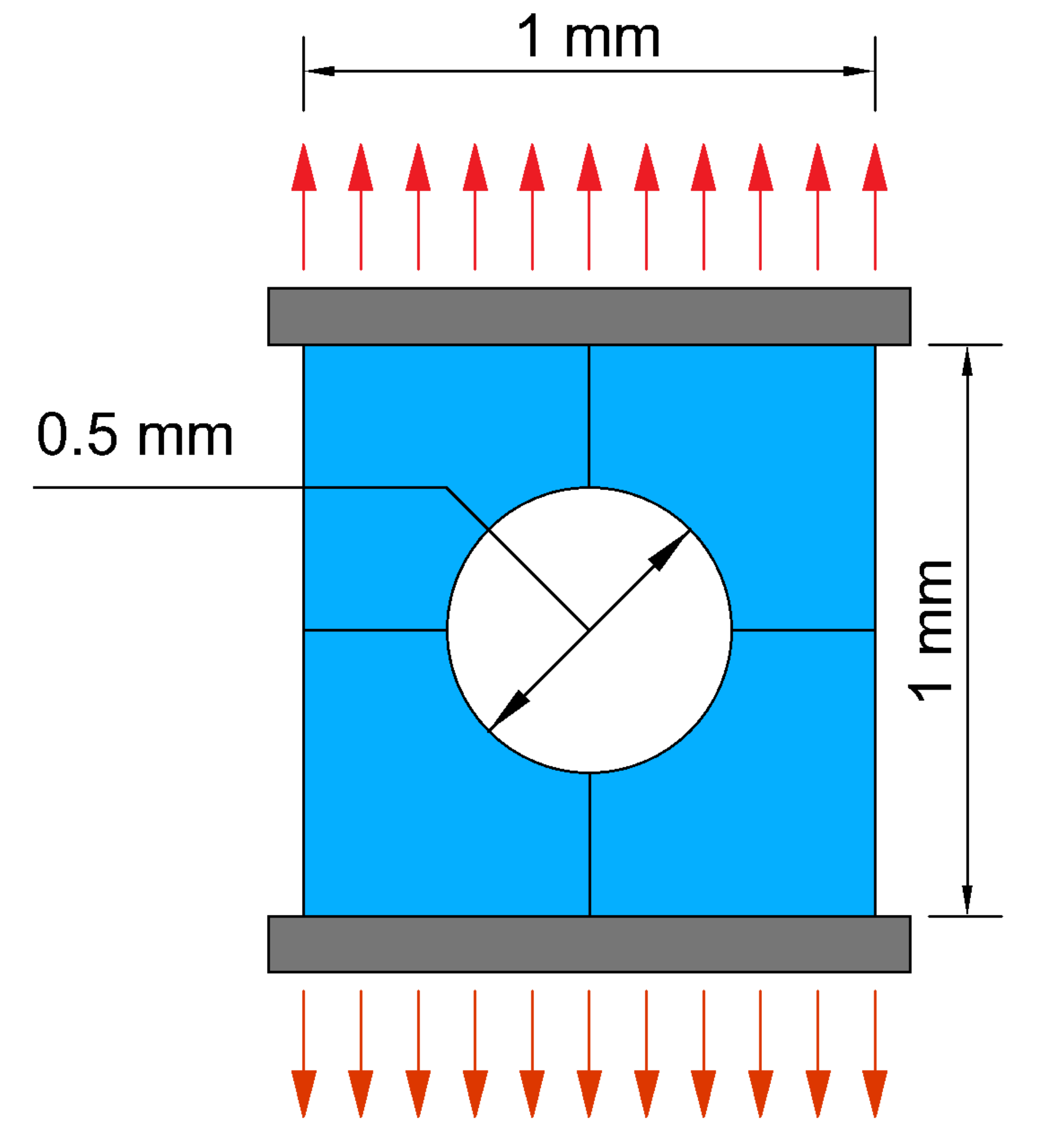}
\vskip 0.1in
\caption{Geometry and boundary conditions for the perforated block under tension: A tensile deformation equal to $300\%$ of the height of the block is applied to the top and bottom edges, which are restrained in the horizontal direction. Due to symmetry in the horizontal and vertical directions, only the upper right quarter of the block is modeled.}
\label{fig:410}
\end{figure}

\begin{figure}
\centering
\includegraphics[width=0.7\textwidth]{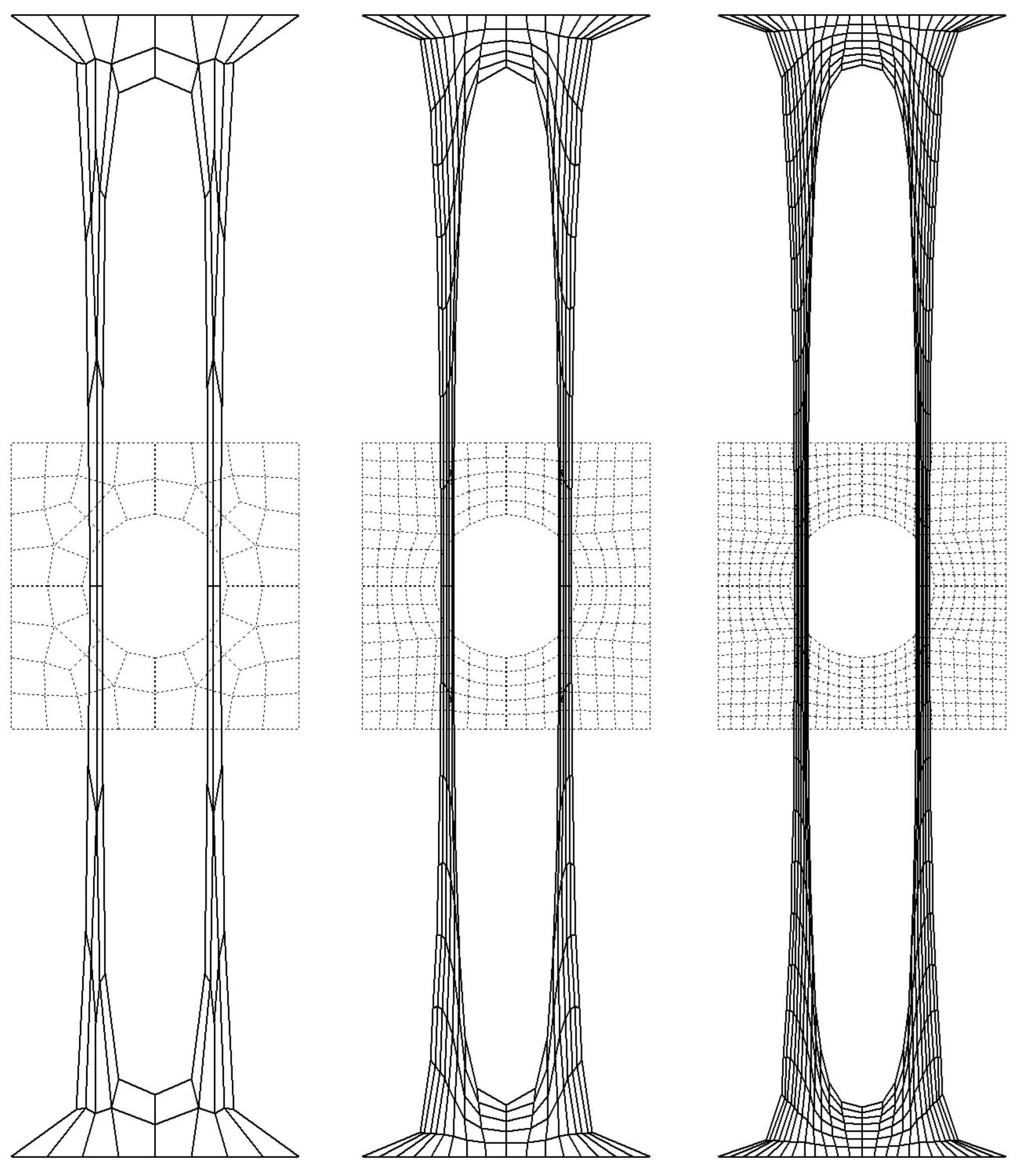}
\vskip 0.1in
\caption{Deformed and undeformed configurations of the perforated block subject to a tensile deformation equal to $300$\% of the height of the block using meshes consisting of (a) $14$, (b) $55$ and (c) $119$ first-order quadrilateral elements developed in this work. The first material model in \eqref{eqn:4.1} with the material parameters $\mu=80.24$ and $\kappa=40,0933.33$~MPa are used in simulations and the top and bottom edges are restrained in the horizontal direction.}
\label{fig:411}
\end{figure}

\begin{figure}
\centering
\includegraphics[width=0.65\textwidth]{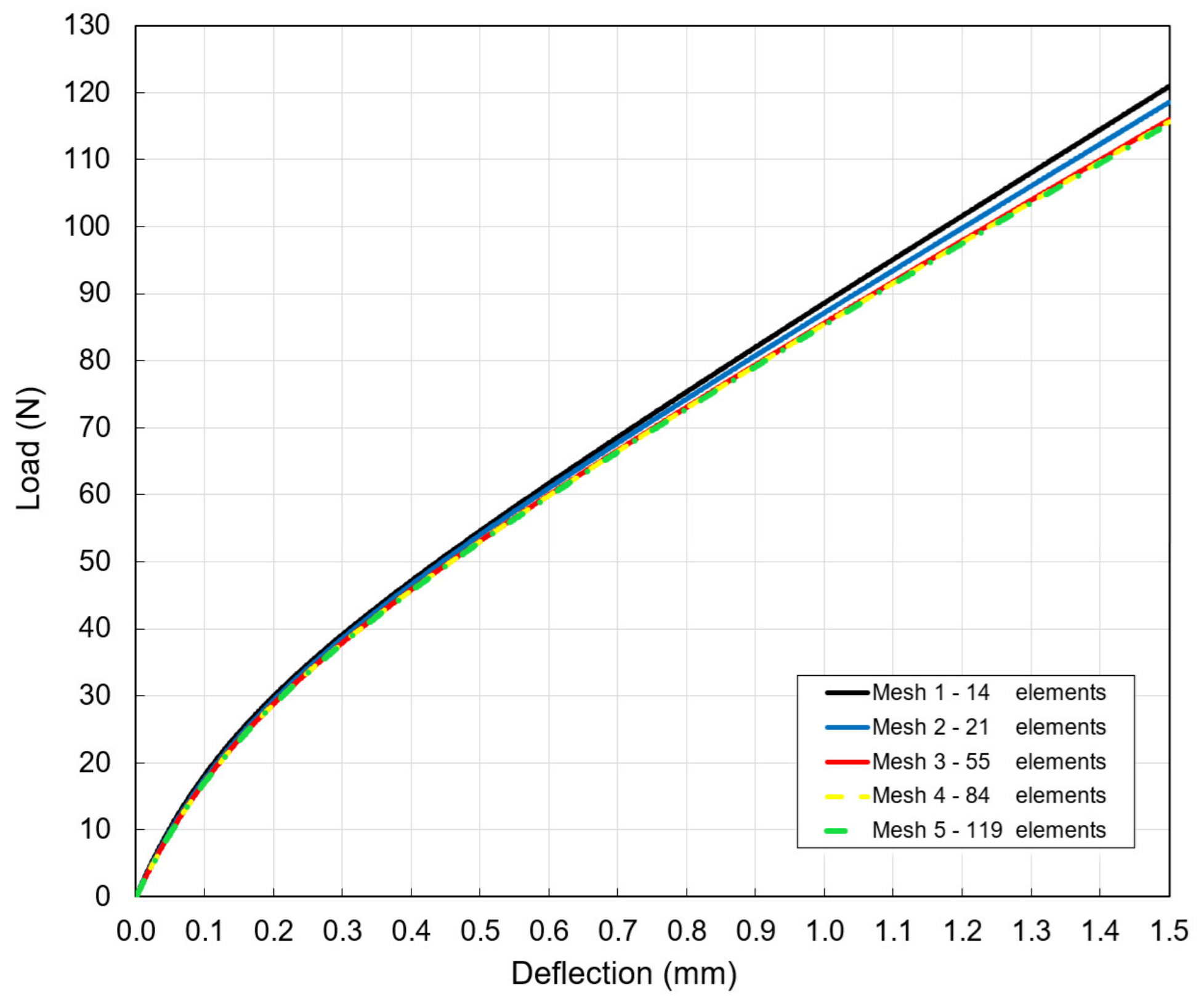}
\vskip 0.1in
\caption{Load-deflection curves for the perforated block subject to a tensile deformation equal to $300$\% of the height of the block using different meshes of the quadrilateral element developed in this work; material parameters are $\mu=80.24$ and $\kappa=40,0933.33$~MPa, and the top and bottom edges are restrained in the horizontal direction.}
\label{fig:412}
\end{figure}

\begin{figure}
\centering
\includegraphics[width=0.6\textwidth]{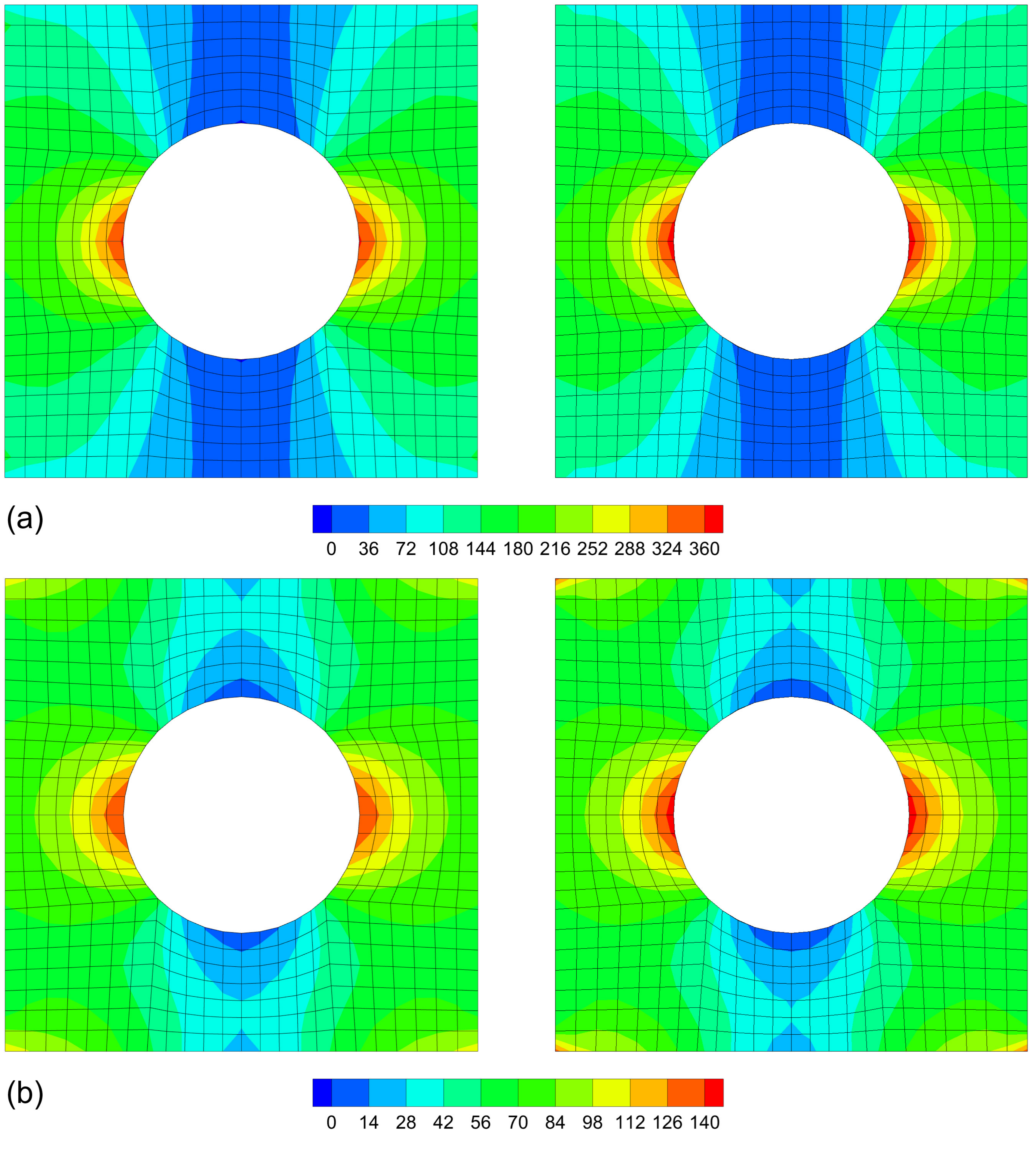}
\vskip 0.1in
\caption{Tension of a perforated block: (a) Contours of the Kirchhoff stress component $\tau_{22}$~(MPa) using the element developed in this work and the second-order mixed U/P finite element. (b) Contours of the pressure $p$~(MPa) using the element developed in this work and the second-order mixed U/P finite element.}
\label{fig:413}
\end{figure}

In this work, a block with a square cross section of $1\times1$ mm and a central hole is subjected to a tensile deformation equal to $300\%$ of its height. The tensile deformation is applied to the top and bottom edges in $1500$ load steps, while both edges are restrained in the horizontal direction. The diameter of the central hole is $0.50$ mm. The geometry and boundary conditions of the problem are shown in Figure~\ref{fig:410}. The simulations are performed using the compressible formulation and the first material model in \eqref{eqn:4.1}, with $\mu=80.24$ and $\kappa=40,093.33$ MPa. Owing to symmetry in both the horizontal and vertical directions, only one quarter of the block is modeled. Several meshes are considered to study the problem. The number of elements, the corresponding number of degrees of freedom, and the final value of the load applied to the top edge for these meshes are summarized in Table~\ref{tab:4.3}. The undeformed and deformed configurations of the block are shown in Figure~\ref{fig:411} for selected meshes. The load-displacement curves resulting from the applied tensile deformation are plotted in Figure~\ref{fig:412} for the meshes listed in Table~\ref{tab:4.3}. It is observed that the curves corresponding to the last three meshes coincide, indicating excellent convergence. The contour plots of the Kirchhoff stress component $\tau_{22}$ and the pressure $p$ in Figures~\ref{fig:413}(a) and \ref{fig:413}(b) are obtained using the incompressible formulation and the second material model in \eqref{eqn:4.1}. In order to compare the results of the incompressible formulation with those of the second-order mixed U/P formulation, the target displacement is limited to $25\%$ of the height of the block. An excellent agreement is observed between the results obtained using the two elements.

\begin{table}
\centering
\caption{The number of elements, the number of degrees of freedom, and the final value of the load applied to the top edge at a target displacement equal to $300\%$ of the height of the block for the meshes of the quadrilateral element developed in this work.}
\label{tab:4.3}
\renewcommand{\arraystretch}{1.5}
\renewcommand{\tabcolsep}{0.2cm}
\begin{tabular}{c c c c}
\hline
Mesh & No. Elements & DOFs   & Load (N)  \\
\hline
1    & 14           & 431    & 120.92    \\
2    & 21           & 629    & 118.56    \\
3    & 55           & 1562   & 115.97    \\
4    & 84           & 2351   & 115.62    \\
5    & 119          & 3295   & 115.39    \\
\hline
\end{tabular}
\end{table}

\clearpage
\section{Conclusions}
\label{sec:5}

In this paper, first-order compatible-strain mixed quadrilateral finite elements were developed for compressible and incompressible nonlinear elasticity. For general quadrilateral elements in the physical space, vector fields that are tangent and normal to an edge in the natural coordinate system do not, in general, remain tangent and normal to the corresponding edge after application of the Piola transformation. To overcome this difficulty, the shape functions used to interpolate the displacement gradient and stress tensor were computed directly in the physical space using numerical integration. As demonstrated by the numerical examples, this feature of the formulation is essential for accommodating highly irregular meshes. This stands in contrast to formulations based on the Piola transformation, which require a certain degree of regularity in the mesh in order to preserve the desired interpolation properties. The stress interpolation and its associated degrees of freedom were designed so that the resulting formulation satisfies the necessary and sufficient conditions for invertibility of the global stiffness matrix.

The present work extends the compatible-strain mixed finite element framework from simplicial to quadrilateral elements for both compressible and incompressible nonlinear elasticity. For incompressible elasticity, the Hu--Washizu functional was augmented by the dilation and pressure fields, while the pressure variable was condensed out at the element level so that no additional global degrees of freedom were introduced. The performance of the proposed elements was assessed through several numerical examples. The load--deflection responses and convergence studies demonstrate excellent convergence properties, particularly for nearly incompressible solids. Accurate solutions were obtained even on relatively coarse meshes. In addition, the proposed quadrilateral elements successfully solved problems that require second-order simplicial CSMFEs while using fewer degrees of freedom. Overall, the proposed quadrilateral elements combine compatible interpolation, numerical stability, and favorable convergence characteristics for the simulation of nonlinear elasticity problems.

Future work will focus on extending the present framework to higher-order quadrilateral elements and three-dimensional hexahedral elements. Another important direction is the development of compatible-strain mixed formulations for finite-strain inelasticity, including elastoplasticity and growth mechanics. It would also be of interest to investigate alternative interpolation spaces and integration strategies that further improve computational efficiency while preserving the compatibility properties of the formulation.

\bibliographystyle{abbrvnat}
\bibliography{ref}

\end{document}